\documentclass[twoside,11pt,reqno]{amsart}
\usepackage{amsmath,amssymb,amscd,mathrsfs,epic,empheq, enumerate}
\usepackage{latexsym,amsthm,amscd}

\makeatletter

\@addtoreset{equation}{section}

\newtheorem{Proposition}{Proposition}[section]
\newtheorem{Lemma}[Proposition]{Lemma}
\newtheorem{Theorem}[Proposition]{Theorem}
\newtheorem{Corollary}[Proposition]{Corollary}

\newtheorem{Remark}[Proposition]{Remark}

\newbox\squ  % box character for ends of proofs
\setbox\squ=\hbox{\vrule width.6pt
   \vbox{\hrule height.6pt width.4em\kern1ex
         \hrule height.6pt}%
   \vrule width.6pt}

\def\phantomsubsection#1{\vspace{2mm}\noindent{\bf #1.}}

\def\deg{\operatorname{deg}}
\def\defect{\operatorname{def}}
\def\cl{{\operatorname{cl}}}
\def\ex{{\operatorname{ex}}}
\def\op{\operatorname{op}}

\def\down{\vee}
\def\up{\wedge}
\def\rep#1{\operatorname{Rep}(#1)}
\def\proj#1{\operatorname{Proj}(#1)}
\def\mod#1{\operatorname{Mod}_{l\!f}(#1)}

\def\C{{\mathbb F}}

\def\Z{{\mathbb Z}}

\def\0{{\bar 0}}
\def\1{{\bar 1}}

\def\hom{{\operatorname{Hom}}}

\def\phi{{\varphi}}

\def\la{{\lambda}}
\def\La{{\Lambda}}
\def\Ga{{\Gamma}}
\def\ga{{\gamma}}
\def\Om{{\Omega}}

\def\De{{\Delta}}
\def\al{{\alpha}}
\def\be{{\beta}}

\newdimen\hoogte    \hoogte=8pt    % hoogte  van hokje
\newdimen\breedte   \breedte=8pt   % breedte van hokje
\newdimen\dikte     \dikte=0.5pt    % dikte lijn

\newenvironment{young}{\begingroup
       \def\vr{\vrule height0.8\hoogte width\dikte depth 0.2\hoogte}
       \def\fbox##1{\vbox{\offinterlineskip
                    \hrule height\dikte
                    \hbox to \breedte{\vr\hfill##1\hfill\vr}
                    \hrule height\dikte}}
       \vbox\bgroup \offinterlineskip \tabskip=-\dikte \lineskip=-\dikte
            \halign\bgroup &\fbox{##\unskip}\unskip  \crcr }
       {\egroup\egroup\endgroup}
\def\diagram#1{\relax\ifmmode\vcenter{\,\begin{young}#1\end{young}\,}\else%
              $\vcenter{\,\begin{young}#1\end{young}\,}$\fi}

\newdimen\Hoogte    \Hoogte=12pt    % hoogte  van hokje
\newdimen\Breedte   \Breedte=12pt   % breedte van hokje
\newdimen\Dikte     \Dikte=0.5pt    % dikte lijn

\newenvironment{Young}{\begingroup
       \def\vr{\vrule height0.8\Hoogte width\Dikte depth 0.2\Hoogte}
       \def\fbox##1{\vbox{\offinterlineskip
                    \hrule height\Dikte
                    \hbox to \Breedte{\vr\hfill##1\hfill\vr}
                    \hrule height\Dikte}}
       \vbox\bgroup \offinterlineskip \tabskip=-\Dikte \lineskip=-\Dikte
            \halign\bgroup &\fbox{##\unskip}\unskip  \crcr }
       {\egroup\egroup\endgroup}
\def\Diagram#1{\relax\ifmmode\vcenter{\,\begin{Young}#1\end{Young}\,}\else%
              $\vcenter{\,\begin{Young}#1\end{Young}\,}$\fi}

\begin{document}

\title[Khovanov's diagram algebra I]{Highest weight categories
arising from Khovanov's diagram algebra I: cellularity}
\author{Jonathan Brundan and Catharina Stroppel}

\address{Department of Mathematics, University of Oregon, Eugene, OR 97403, USA}
\email{brundan@uoregon.edu}
\address{Department of Mathematics, University of Bonn, 53115 Bonn, Germany}
\email{stroppel@math.uni-bonn.de}

\thanks{First author supported in part by NSF grant no. DMS-0654147.}
\thanks{Second author supported by the NSF and the Minerva Research Foundation DMS-0635607.}

\begin{abstract}
This is the first of four articles
studying some slight generalisations $H^n_m$
of Khovanov's diagram algebra,
as well as  quasi-hereditary covers
$K^n_m$ of these algebras in the sense of Rouquier, and certain
infinite dimensional limiting versions $K^\infty_m$,
$K^{\pm \infty}_m$ and $K^\infty_\infty$.
In this article we prove that $H^n_m$ is a cellular symmetric algebra
and that $K^n_m$ is a cellular quasi-hereditary
algebra. In subsequent articles, we
relate $H^n_m, K^n_m$ and $K^\infty_m$
to level two blocks of
degenerate cyclotomic Hecke algebras,
parabolic category $\mathcal O$ and the general linear supergroup,
respectively.
\end{abstract}
\maketitle

\tableofcontents

\section{Introduction}\label{sintroduction}

Quasi-hereditary algebras (and highest weight categories) play
an important role in representation theory thanks to their nice homological properties. They are finite dimensional algebras with distinguished
collections of projective indecomposable, standard and irreducible modules, all
of which are labeled by the same indexing set.
Module categories over quasi-hereditary algebras are good candidates for
categorifications, since each of these three collections of modules defines a
distinguished basis for the Grothendieck group, and the transition matrices
between these bases are triangular.
Quasi-hereditariness implies finite global dimension, and
so excludes many interesting algebras like
group algebras or Hecke algebras. In \cite{GL}
Graham
and Lehrer introduced the weaker notion of a
{cellular algebra} which again comes along with three distinguished classes of
modules.
Cellular algebras often have a more combinatorial nature, and
include many diagram algebras such as Temperley-Lieb algebras and
Brauer algebras which appear in connection with knot invariants.

Khovanov's categorification of the Jones polynomial is based on
certain diagrammatically defined algebras
denoted here by $H_n^n$, whose
multiplication is defined via some 2-dimensional TQFT structure;
see \cite{K}, \cite{K2}.
This article is the first of a series of {four} articles
studying some slight
generalisations $H^n_m$ of Khovanov's diagram
algebra, and their quasi-hereditary covers $K^n_m$ which
were introduced already in \cite{ChK}, \cite{S}.
All of these are positively graded finite dimensional algebras;
since
$H^n_m \cong H^m_n$ and $K^n_m \cong K^m_n$ we can assume
that $m \leq n$.
By taking
direct limits, we also define and study some interesting
positively graded infinite dimensional
algebras denoted
$K^\infty_m$, $K^{\pm \infty}_m$ and $K^\infty_\infty$.

Here is a very brief summary of the results to be proved about these
algebras in this article and its sequels.

\vspace{2mm}

In Part I (this article) we show
that the algebras $H^n_m$
are cellular symmetric algebras and
the algebras $K^n_m$ are cellular quasi-hereditary algebras.
We also compute explicitly the {$q$-decomposition matrices}
that describe the transition between the various distinguished bases
for the graded Grothendieck groups.
Analogous results are obtained for the limits
$K^\infty_m$ and $K^{\pm \infty}_m$ too, though because they
are locally unital but not unital algebras
we must say instead that their module categories are highest weight categories.
Amusingly, it happens that $K^\infty_m$
and $K^{+\infty}_m$ are also symmetric algebras,
a combination which would be impossible in the finite dimensional setting; see also the recent work \cite{CT} for
another class of locally unital algebras
with similar properties.

In Part II we study the appropriate analogues of
Khovanov's diagrammatically defined
{``projective functors''} for our various families of algebras.
In Khovanov's work, these functors play the central role in the
categorification of the Jones polynomial.
They play an equally central but very different role in our work.
For example, we use them to give an elementary proof that the algebras
$K^n_m, K_m^\infty, K_m^{\pm \infty}$ and $K^\infty_\infty$ are Koszul and compute the
associated Kazhdan-Lusztig polynomials explicitly;
these polynomials are
the entries of the inverse of the $q$-decomposition matrix
mentioned in the previous paragraph.
Projective functors also allow us to prove the double centraliser property
that justifies our earlier
statement that $K^n_m$ is a quasi-hereditary cover
of $H^n_m$.

In Part III, we give a new proof of a known result
relating the algebra $K^n_m$ to the blocks of a certain
category $\mathcal O(m,n)$, namely, the
parabolic analogue of the Bernstein-Gelfand-Gelfand
category $\mathcal O$ associated to the
subalgebra $\mathfrak{gl}_m(\mathbb{C}) \oplus \mathfrak{gl}_n(\mathbb{C})$
in $\mathfrak{gl}_{m+n}(\mathbb{C})$ (we consider only integral central characters).
In \cite{S},
the second author showed
that the algebra $K^n_m$
is isomorphic to the endomorphism algebra of a projective
generator for the category of perverse sheaves on the Grassmannian
of $m$-dimensional
subspaces in an $(m+n)$-dimensional complex vector space.
This was proved by identifying the diagrammatically defined algebra
$K^n_m$ with a very different presentation
of this endomorphism algebra due to Braden \cite{Br}.
Combined with the Beilinson-Bernstein localisation theorem and
the Riemann-Hilbert correspondence,
this means that $K^n_m$ is also the endomorphism algebra of a projective
generator for a regular block of $\mathcal O(m,n)$.
Our new approach to the proof of this
statement bypasses
geometry entirely: it is based
instead on the generalised
Schur-Weyl duality from \cite{BK} which
connects $\mathcal O(m,n)$
to degenerate cyclotomic Hecke algebras of level two.
We actually prove a more general result
giving an equivalence between the entire category
$\mathcal O(m,n)$
and a category
introduced
in \cite{HK}, \cite{Ch}.
As a by-product
we recover an old result
of Enright and Shelton
formulated in a special case
in \cite[Proposition 11.2]{ES}: singular blocks of $\mathcal O(m,n)$
are equivalent to regular blocks
for smaller $m$ and $n$.
Moreover, we show that
the generalised Khovanov algebras $H^n_m$ are Morita equivalent
to the blocks of
degenerate cyclotomic Hecke algebras of level two,
and prove
a very recent conjecture of Khovanov and Lauda from \cite[$\S$3.4]{KLa}
in the case of level two weights in finite type A.

In Part IV, we
relate the limiting versions $K^\infty_r$ to blocks of the general
linear supergroup
$GL(m|n)$, a setting where a geometric approach via some sort of localisation
theorem is missing.
More precisely, we prove that the (locally finite) endomorphism algebra of a projective
generator for a block of $GL(m|n)$
of atypicality $r$
is isomorphic to the
algebra $K^\infty_r$.
Combined with results from Part II,
this implies that blocks of $GL(m|n)$
are Koszul.
It also shows that all blocks of the same atypicality
are equivalent, which is an old unpublished result of Serganova.
The strategy for the proof in Part IV
is similar in spirit to the argument in Part III,
being based on a new ``super''
Schur-Weyl duality
relating the general linear supergroup to the same
degenerate cyclotomic Hecke algebras of level two
mentioned before.
Finally, we mention that by combining the results of Parts III and IV, one can deduce a proof of a recent conjecture of
Cheng, Wang and Zhang \cite{CWZ}. In fact this conjecture was the
original motivation leading us to the limiting versions
of Khovanov's algebra.

\vspace{2mm}

\iffalse
We remark that
some
signed versions
$\widetilde{H}_m^n$ and $\widetilde{K}_m^n$
of the algebras $H_m^n$ and $K_m^n$ have recently been introduced
in \cite{SW}.
They
appear naturally in the geometric incarnations of Khovanov homology via coherent sheaves \cite{CK} or symplectic geometry
related to Slodowy slices \cite{SS}.
As explained in \cite[$\S$5.2]{SW}, the construction
of the signed versions differs from
Khovanov's construction
in that the underlying 2-dimensional TQFT structure is enhanced
by working inside the category of disoriented cobordisms in the sense of
\cite{CMW}.
The methods of this article (and also those of
Part II) carry over with remarkably few changes to yield similar
results for these signed versions.
In particular, the combinatorics underlying the representation theory
is essentially the same in both cases.
\fi

Since the first version of this article became available, there have been
several subsequent developments, some of which we mention briefly here.
For an interpretation of Khovanov's algebra in terms of the geometry of the
Springer fibre we refer to \cite{SW}. In that work,
weights are interpreted as fixed points of a torus action on a
Springer fibre,
the oriented circle diagrams defined below correspond to fixed points lying inside
the intersection $C_1\cap C_2$ of two irreducible
components hence
index a basis of $H^*(C_1\cap C_2)$, and Khovanov's
combinatorial multiplication rule
arises from some convolution product.
On the category $\mathcal O$ side (Paper III), our diagrammatics has
recently been extended to the Hermitian symmetric pair $(A_{n-1}, D_{n})$
by Lejcyk \cite{Lej}.
Intriguingly the diagrams in \cite{Lej} are very similar to ones
appearing in another recent paper by Cox-De
Visscher \cite{CD} on combinatorial representation theory of
(walled) Brauer algebras.
On the supergroup side (Paper IV), a more direct
proof of the super duality conjecture of Cheng, Wang and Zhang has
now been found by Cheng and Lam \cite{CL}.
Also our diagrammatics for $GL(m|n)$ has been
extended to cover blocks of the orthosymplectic supergroup $OSp(m|2n)$
by Gruson and Serganova \cite{CG}.

\vspace{2mm}
We devote the remainder of this introduction to giving an informal
guide to the construction of the algebras $H^n_m$, $K^n_m$ and their
limiting versions. We hope the baby examples sketched
here make the precise definitions explained in the main body of the
article more accessible.
As a vector space, Khovanov's algebra $H^n_n$ is usually identified with
$$
\bigoplus_{ab} \left[\mathbb{C}[x] / (x^2)\right]^{\otimes  (\text{the number of circles in $ab$})}
$$
where the direct sum is over all closed circle diagrams $ab$
obtained by gluing a cup diagram $a$
with $n$ cups {under} a cap diagram $b$ with $n$
caps;
see $\S$2 for all this language.
The summand of $H^n_n$ corresponding to the circle diagram $ab$ has a
distinguished basis consisting
of tensors in $1$ and $x$.
In the literature,
this distinguished basis is often parametrised by decorating each
circle in $ab$ by an additional {type}, $1$ or $x$.
In this article, we instead encode circles of types $1$ or $x$
by {orienting} them so that they are anti-clockwise or clockwise, respectively.
For example, we represent the distinguished
basis $\{1,x\}$ of $H^1_1 = \mathbb{C}[x] / (x^2)$ by the
closed oriented circle diagrams
$$
\begin{picture}(-20,30)
\put(18.2,11.6){$x \:\,\leftrightarrow$}
\put(49.2,9.4){$\scriptstyle\up$}
\put(72.2,14.1){$\scriptstyle\down$}
\put(63.5,14){\oval(23,23)[b]}
\put(63.5,14){\oval(23,23)[t]}
\put(75,14){\line(-1,0){23}}

\put(-109.8,11.6){$1 \:\,\leftrightarrow$}
\put(-63.5,14){\oval(23,23)[b]}
\put(-63.5,14){\oval(23,23)[t]}
\put(-75,14){\line(1,0){23}}
\put(-54.8,9.4){$\scriptstyle\up$}
\put(-77.8,14.1){$\scriptstyle\down$}
\end{picture}
$$
In fact $H^1_1$ equals $\mathbb{C}[x] / (x^2)$
{\em as an algebra}, making it a deceptively simple example.
Khovanov's definition of the multiplication
on $H^n_n$ in general involves a 2-dimensional TQFT
arising from $\mathbb{C}[x] / (x^2)$.
In this article, we have reformulated Khovanov's definition
in a slightly different way.
Here is an example giving the product of two basis vectors
in $H^2_2$:
$$
\begin{picture}(8,50)
\put(-133.5,24){\oval(23,23)[b]}
\put(-133.5,24){\oval(69,46)[b]}
\put(-110.5,24){\oval(23,23)[t]}
\put(-156.5,24){\oval(23,23)[t]}
\put(-168,24){\line(1,0){69}}
\put(-124.8,24.1){$\scriptstyle\down$}
\put(-101.8,19.4){$\scriptstyle\up$}
\put(-147.8,19.4){$\scriptstyle\up$}
\put(-170.8,24.1){$\scriptstyle\down$}
\put(-93,21.5){$\times$}
\put(-43.5,24){\oval(23,23)[t]}
\put(-43.5,24){\oval(69,46)[t]}
\put(-20.5,24){\oval(23,23)[b]}
\put(-66.5,24){\oval(23,23)[b]}
\put(-78,24){\line(1,0){69}}
\put(-34.8,24.1){$\scriptstyle\down$}
\put(-11.8,19.4){$\scriptstyle\up$}
\put(-57.8,19.4){$\scriptstyle\up$}
\put(-80.8,24.1){$\scriptstyle\down$}
\put(-0.5,21.5){$=$}
\put(51.5,24){\oval(23,23)[t]}
\put(51.5,24){\oval(69,46)[t]}
\put(51.5,24){\oval(23,23)[b]}
\put(51.5,24){\oval(69,46)[b]}
\put(86,24){\line(-1,0){69}}
\put(37.2,19.4){$\scriptstyle\up$}
\put(14.2,24.1){$\scriptstyle\down$}
\put(60.2,24.1){$\scriptstyle\down$}
\put(83.2,19.4){$\scriptstyle\up$}
\put(92,21.5){$+$}
\put(141.5,24){\oval(23,23)[t]}
\put(141.5,24){\oval(69,46)[t]}
\put(141.5,24){\oval(23,23)[b]}
\put(141.5,24){\oval(69,46)[b]}
\put(176,24){\line(-1,0){69}}
\put(127.2,24.1){$\scriptstyle\down$}
\put(104.2,19.4){$\scriptstyle\up$}
\put(150.2,19.4){$\scriptstyle\up$}
\put(173.2,24.1){$\scriptstyle\down$}
\end{picture}
$$
Roughly, to calculate this product in our setup,
the idea is to
draw the first diagram under the second and then perform
two ``surgery procedures'' following rules explained in detail in $\S$3.

Like in the above pictures, a choice of orientation on
a circle diagram
gives rise to a {``weight''}, that is, a sequence
of $\down$'s and $\up$'s on
the horizontal ``number line'' that separates cups from caps.
So it makes sense to denote oriented circle diagrams
by triples $a \la b$, where $a$ is the cup diagram at the bottom,
$\la$ is the weight in the middle encoding the orientation,
and $b$ is the cap diagram at the top.
The weights arising from $H^n_n$ have $n$ $\down$'s
and $n$ $\up$'s;
for example, the set of all weights for $H^2_2$ is
$\La :=
\{{{\down\!\down\!\up\up}},\:
{{\down\!\up\!\down\up}},\:
{{\up\!\down\!\down\up}},\:
{{\down\!\up\!\up\down}},\:
{{\up\!\down\!\up\down}},\:
{{\up\!\up\!\down\down}}
\}.
$
There is a partial order $\leq$ on $\La$
generated by swapping $\down \up$ pairs.
In the $H_2^2$ example displayed in the previous paragraph,
the weight ${\up\!\down\!\up\down}$
on the right hand side
is strictly greater than the weight
${\down\!\up\!\down\up}$ on the left hand side.
This is a general phenomenon: it is {\em always}
the case that the product
$(a \la b) (c \mu d)$ of two basis vectors
is a linear combination $(a \nu d)$ of basis vectors
with $\la \leq \nu \geq \mu$.
This leads to the cellular structure on $H^n_n$;
see $\S$3.

For $m < n$, the generalised Khovanov algebra $H^n_m$ has a very
similar cellular basis parametrised by (no longer closed)
oriented circle
diagrams $a \la b$.
The cup diagram $a$ at the bottom
now has $m$ cups and $(n-m)$ rays going down to infinity,
the weight $\la$ in the middle
encoding the orientation has $m$ $\down$'s and $n$ $\up$'s,
and the cap diagram $b$ at the top
has $m$ caps and $(n-m)$ rays going up to infinity.
For example, $H^2_1$ has basis
parametrised by the following six diagrams:
$$
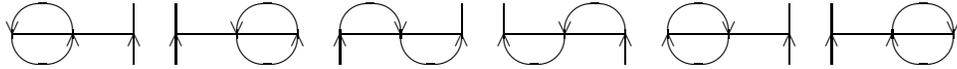
\begin{picture}(62,28)
\put(9.2,9.4){$\scriptstyle\up$}
\put(32.2,9.4){$\scriptstyle\up$}
\put(-13.8,14.1){$\scriptstyle\down$}
\put(.5,14){\oval(23,23)[b]}
\put(.5,14){\oval(23,23)[t]}
\put(35,2.5){\line(0,1){23}}
\put(35,14){\line(-1,0){46}}
\end{picture}
\begin{picture}(62,28)
\put(9.2,14.1){$\scriptstyle\down$}
\put(32.2,9.4){$\scriptstyle\up$}
\put(-13.8,9.4){$\scriptstyle\up$}
\put(23.5,14){\oval(23,23)[b]}
\put(23.5,14){\oval(23,23)[t]}
\put(-11,2.5){\line(0,1){23}}
\put(35,14){\line(-1,0){46}}
\end{picture}
\begin{picture}(62,28)
\put(9.2,14.1){$\scriptstyle\down$}
\put(32.2,9.4){$\scriptstyle\up$}
\put(-13.8,9.4){$\scriptstyle\up$}
\put(23.5,14){\oval(23,23)[b]}
\put(.5,14){\oval(23,23)[t]}
\put(-11,14){\line(0,-1){11.5}}
\put(35,14){\line(0,1){11.5}}
\put(35,14){\line(-1,0){46}}
\end{picture}
\begin{picture}(62,28)
\put(9.2,14.1){$\scriptstyle\down$}
\put(32.2,9.4){$\scriptstyle\up$}
\put(-13.8,9.4){$\scriptstyle\up$}
\put(23.5,14){\oval(23,23)[t]}
\put(.5,14){\oval(23,23)[b]}
\put(-11,14){\line(0,1){11.5}}
\put(35,14){\line(0,-1){11.5}}
\put(35,14){\line(-1,0){46}}
\end{picture}
\begin{picture}(62,28)
\put(9.2,14.1){$\scriptstyle\down$}
\put(32.2,9.4){$\scriptstyle\up$}
\put(-13.8,9.4){$\scriptstyle\up$}
\put(.5,14){\oval(23,23)[b]}
\put(.5,14){\oval(23,23)[t]}
\put(35,2.5){\line(0,1){23}}
\put(35,14){\line(-1,0){46}}
\end{picture}
\begin{picture}(30,26)
\put(9.2,9.4){$\scriptstyle\up$}
\put(32.2,14.1){$\scriptstyle\down$}
\put(-13.8,9.4){$\scriptstyle\up$}
\put(23.5,14){\oval(23,23)[b]}
\put(23.5,14){\oval(23,23)[t]}
\put(-11,2.5){\line(0,1){23}}
\put(35,14){\line(-1,0){46}}
\end{picture}
$$
The multiplication on $H^m_n$ can be defined
by a slightly
generalised surgery procedure; see $\S$6 for details.
We have not yet mentioned the natural $\Z$-grading on all our algebras.
This can be defined simply by declaring that the degree of a basis vector
is equal to the total number of clockwise cups and caps in its diagram.
The above basis vectors for $H^2_1$ are of degrees
$0,0,1,1,2, 2$. As this example suggests, all $H^n_m$ are symmetric algebras;
see $\S$6.

Now we can outline the construction of the quasi-hereditary covers $K^n_m$
of $H^n_m$.
These also have bases parametrised
by oriented circle diagrams $a \la b$,
and the weights $\la$ decorating the number line
again have $m$ $\down$'s and $n$ $\up$'s.
However we now allow cup and cap diagrams $a$ and $b$ with fewer
cups and caps than required before, compensating by adding additional
rays to infinity. There is also one new
constraint on orientations: we do not allow
two rays to infinity to be oriented
$\down \up$ in that order from left to right.
\begin{table}\label{t1}
\title{Table 1: Multiplication in $K^1_1$}
$$
\begin{array}{l|lllll}
\begin{picture}(30,28)
\put(8,10){$K^1_1$}
\end{picture}
&
\begin{picture}(30,28)
\put(9.2,9.4){$\scriptstyle\up$}
\put(32.2,14.1){$\scriptstyle\down$}
\put(35,2.5){\line(0,1){23}}
\put(12,2.5){\line(0,1){23}}
\put(35,14){\line(-1,0){23}}
\end{picture}
&
\begin{picture}(30,28)
\put(9.2,9.4){$\scriptstyle\up$}
\put(32.2,14.1){$\scriptstyle\down$}
\put(35,14){\line(0,-1){11.5}}
\put(12,14){\line(0,-1){11.5}}
\put(23.5,14){\oval(23,23)[t]}
\put(35,14){\line(-1,0){23}}
\end{picture}
&
\begin{picture}(30,28)
\put(9.2,9.4){$\scriptstyle\up$}
\put(32.2,14.1){$\scriptstyle\down$}
\put(35,14){\line(0,1){11.5}}
\put(12,14){\line(0,1){11.5}}
\put(23.5,14){\oval(23,23)[b]}
\put(35,14){\line(-1,0){23}}
\end{picture}
&
\begin{picture}(30,28)
\put(9.2,14.1){$\scriptstyle\down$}
\put(32.2,9.4){$\scriptstyle\up$}
\put(23.5,14){\oval(23,23)[b]}
\put(23.5,14){\oval(23,23)[t]}
\put(35,14){\line(-1,0){23}}
\end{picture}
&
\begin{picture}(40,28)
\put(23.5,14){\oval(23,23)[b]}
\put(23.5,14){\oval(23,23)[t]}
\put(35,14){\line(-1,0){23}}
\put(9.2,9.4){$\scriptstyle\up$}
\put(32.2,14.1){$\scriptstyle\down$}
\end{picture}
%endfirstline
\\
\hline
\\
\begin{picture}(30,28)
\put(-0.8,9.4){$\scriptstyle\up$}
\put(22.2,14.1){$\scriptstyle\down$}
\put(25,2.5){\line(0,1){23}}
\put(2,2.5){\line(0,1){23}}
\put(25,14){\line(-1,0){23}}
\end{picture}
&
\begin{picture}(30,28)
\put(9.2,9.4){$\scriptstyle\up$}
\put(32.2,14.1){$\scriptstyle\down$}
\put(35,2.5){\line(0,1){23}}
\put(12,2.5){\line(0,1){23}}
\put(35,14){\line(-1,0){23}}
\end{picture}
&
\begin{picture}(30,28)
\put(9.2,9.4){$\scriptstyle\up$}
\put(32.2,14.1){$\scriptstyle\down$}
\put(35,14){\line(0,-1){11.5}}
\put(12,14){\line(0,-1){11.5}}
\put(23.5,14){\oval(23,23)[t]}
\put(35,14){\line(-1,0){23}}
\end{picture}
&
\begin{picture}(30,28)
\put(20.7,9.4){$0$}
\end{picture}
&
\begin{picture}(30,28)
\put(20.7,9.4){$0$}
\end{picture}
&
\begin{picture}(40,28)
\put(20.7,9.4){$0$}
\end{picture}
\\
%endsecondline
\begin{picture}(28,28)
\put(-0.8,9.4){$\scriptstyle\up$}
\put(22.2,14.1){$\scriptstyle\down$}
\put(25,14){\line(0,-1){11.5}}
\put(2,14){\line(0,-1){11.5}}
\put(13.5,14){\oval(23,23)[t]}
\put(25,14){\line(-1,0){23}}
\end{picture}
&
\begin{picture}(30,28)
\put(20.7,9.4){$0$}
\end{picture}
&
\begin{picture}(30,28)
\put(20.7,9.4){$0$}
\end{picture}
&
\begin{picture}(30,28)
\put(20.7,9.4){$0$}
\end{picture}
&
\begin{picture}(30,28)
\put(9.2,9.4){$\scriptstyle\up$}
\put(32.2,14.1){$\scriptstyle\down$}
\put(35,14){\line(0,-1){11.5}}
\put(12,14){\line(0,-1){11.5}}
\put(23.5,14){\oval(23,23)[t]}
\put(35,14){\line(-1,0){23}}
\end{picture}
&
\begin{picture}(40,28)
\put(20.7,9.4){$0$}
\end{picture}
\\
%endthirdline
\begin{picture}(28,28)
\put(-0.8,9.4){$\scriptstyle\up$}
\put(22.2,14.1){$\scriptstyle\down$}
\put(25,14){\line(0,1){11.5}}
\put(2,14){\line(0,1){11.5}}
\put(13.5,14){\oval(23,23)[b]}
\put(25,14){\line(-1,0){23}}
\end{picture}
&
\begin{picture}(30,28)
\put(9.2,9.4){$\scriptstyle\up$}
\put(32.2,14.1){$\scriptstyle\down$}
\put(35,14){\line(0,1){11.5}}
\put(12,14){\line(0,1){11.5}}
\put(23.5,14){\oval(23,23)[b]}
\put(35,14){\line(-1,0){23}}
\end{picture}
&
\begin{picture}(30,28)
\put(9.2,9.4){$\scriptstyle\up$}
\put(32.2,14.1){$\scriptstyle\down$}
\put(23.5,14){\oval(23,23)[b]}
\put(23.5,14){\oval(23,23)[t]}
\put(35,14){\line(-1,0){23}}
\end{picture}
&
\begin{picture}(30,28)
\put(20.7,9.4){$0$}
\end{picture}
&
\begin{picture}(30,28)
\put(20.7,9.4){$0$}
\end{picture}
&
\begin{picture}(40,28)
\put(20.7,9.4){$0$}
\end{picture}
\\
%endfourthline
\\
\begin{picture}(30,20)
\put(-0.8,14.1){$\scriptstyle\down$}
\put(22.2,9.4){$\scriptstyle\up$}
\put(13.5,14){\oval(23,23)[b]}
\put(13.5,14){\oval(23,23)[t]}
\put(25,14){\line(-1,0){23}}
\end{picture}
&
\begin{picture}(30,20)
\put(20.7,9.4){$0$}
\end{picture}
&
\begin{picture}(30,20)
\put(20.7,9.4){$0$}
\end{picture}
&
\begin{picture}(30,20)
\put(9.2,9.4){$\scriptstyle\up$}
\put(32.2,14.1){$\scriptstyle\down$}
\put(35,14){\line(0,1){11.5}}
\put(12,14){\line(0,1){11.5}}
\put(23.5,14){\oval(23,23)[b]}
\put(35,14){\line(-1,0){23}}
\end{picture}
&
\begin{picture}(30,20)
\put(9.2,14.1){$\scriptstyle\down$}
\put(32.2,9.4){$\scriptstyle\up$}
\put(23.5,14){\oval(23,23)[b]}
\put(23.5,14){\oval(23,23)[t]}
\put(35,14){\line(-1,0){23}}
\end{picture}
&
\begin{picture}(40,20)
\put(23.5,14){\oval(23,23)[b]}
\put(23.5,14){\oval(23,23)[t]}
\put(35,14){\line(-1,0){23}}
\put(9.2,9.4){$\scriptstyle\up$}
\put(32.2,14.1){$\scriptstyle\down$}
\end{picture}
\\
\begin{picture}(28,28)
\put(-0.8,9.4){$\scriptstyle\up$}
\put(22.2,14.1){$\scriptstyle\down$}
\put(13.5,14){\oval(23,23)[b]}
\put(13.5,14){\oval(23,23)[t]}
\put(25,14){\line(-1,0){23}}
\end{picture}
&
\begin{picture}(30,28)
\put(20.7,9.4){$0$}
\end{picture}
&
\begin{picture}(30,28)
\put(20.7,9.4){$0$}
\end{picture}
&
\begin{picture}(30,28)
\put(20.7,9.4){$0$}
\end{picture}
&
\begin{picture}(30,28)
\put(9.2,9.4){$\scriptstyle\up$}
\put(32.2,14.1){$\scriptstyle\down$}
\put(23.5,14){\oval(23,23)[b]}
\put(23.5,14){\oval(23,23)[t]}
\put(35,14){\line(-1,0){23}}
\end{picture}
&
\begin{picture}(40,28)
\put(20.7,9.4){$0$}
\end{picture}
\end{array}
$$
\end{table}
The multiplication on $K^n_m$ is somewhat harder to define than on $H^n_n$; the basic idea going back to Khovanov and Braden
is to construct $K^n_m$ as a quotient of $H^{m+n}_{m+n}$. There is also
a direct way to compute multiplication in $K^m_n$ in terms of the
generalised surgery procedure which is much quicker in practise,
though from a theoretical point of view the quotient
construction is better.
We give an example illustrating the quotient construction in the next paragraph,
referring the reader to $\S$4 for details and a proof of cellularity
and $\S$6 for the quicker way to compute multiplication.
In $\S$5, we show $K^n_m$ is a quasi-hereditary algebra
and construct its projective indecomposable,
standard and irreducible modules explicitly; all three
families are parametrised by the same underlying set $\La$ of weights.

The algebra $K^1_1$ is a certain five dimensional algebra that is
well known to category $\mathcal O$ experts; see Table~1
for its full multiplication table in terms of diagrams. The entry
$$
\begin{picture}(45,28)
\put(9.2,9.4){$\scriptstyle\up$}
\put(32.2,14.1){$\scriptstyle\down$}
\put(35,14){\line(0,1){11.5}}
\put(12,14){\line(0,1){11.5}}
\put(23.5,14){\oval(23,23)[b]}
\put(35,14){\line(-1,0){23}}
\put(41,12){$\times$}
\end{picture}
\begin{picture}(45,28)
\put(9.2,9.4){$\scriptstyle\up$}
\put(32.2,14.1){$\scriptstyle\down$}
\put(35,14){\line(0,-1){11.5}}
\put(12,14){\line(0,-1){11.5}}
\put(23.5,14){\oval(23,23)[t]}
\put(35,14){\line(-1,0){23}}
\put(41,12){$=$}
\end{picture}
\begin{picture}(40,28)
\put(23.5,14){\oval(23,23)[b]}
\put(23.5,14){\oval(23,23)[t]}
\put(35,14){\line(-1,0){23}}
\put(9.2,9.4){$\scriptstyle\up$}
\put(32.2,14.1){$\scriptstyle\down$}
\end{picture}
$$
of this table is computed as follows.
First ``close'' the diagrams to be multiplied
by adding one extra vertex labeled $\down$ to the left end
and one labeled $\up$ to the right end
of the number line,
then
adding some
new anti-clockwise cups and caps
connecting the rays into closed circles through these new
vertices:
$$
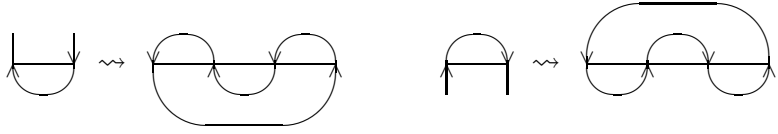
\begin{picture}(76,50)
\put(-33.5,24){\oval(23,23)[b]}
\put(-45,24){\line(0,1){11.5}}
\put(-22,24){\line(0,1){11.5}}
\put(-45,24){\line(1,0){23}}
\put(-24.8,24.1){$\scriptstyle\down$}
\put(-47.8,19.4){$\scriptstyle\up$}
\put(-13,22){$\rightsquigarrow$}
\end{picture}
\begin{picture}(88,50)
\put(-33.5,24){\oval(23,23)[b]}
\put(-33.5,24){\oval(69,46)[b]}
\put(-10.5,24){\oval(23,23)[t]}
\put(-56.5,24){\oval(23,23)[t]}
%\dashline{2.3}(-68,24)(-45,24)
\put(-68,24){\line(1,0){69}}
%\dashline{2.3}(-22,24)(1,24)
\put(-24.8,24.1){$\scriptstyle\down$}
\put(-1.8,19.4){$\scriptstyle\up$}
\put(-47.8,19.4){$\scriptstyle\up$}
\put(-70.8,24.1){$\scriptstyle\down$}
\end{picture}
\begin{picture}(76,50)
\put(-33.5,24){\oval(23,23)[t]}
\put(-45,24){\line(0,-1){11.5}}
\put(-22,24){\line(0,-1){11.5}}
\put(-45,24){\line(1,0){23}}
\put(-24.8,24.1){$\scriptstyle\down$}
\put(-47.8,19.4){$\scriptstyle\up$}
\put(-13,22){$\rightsquigarrow$}
\end{picture}
\begin{picture}(-30,50)
\put(-33.5,24){\oval(23,23)[t]}
\put(-33.5,24){\oval(69,46)[t]}
\put(-10.5,24){\oval(23,23)[b]}
\put(-56.5,24){\oval(23,23)[b]}
%\dashline{2.3}(-68,24)(-45,24)
\put(-68,24){\line(1,0){69}}
%\dashline{2.3}(-22,24)(1,24)
\put(-24.8,24.1){$\scriptstyle\down$}
\put(-1.8,19.4){$\scriptstyle\up$}
\put(-47.8,19.4){$\scriptstyle\up$}
\put(-70.8,24.1){$\scriptstyle\down$}
\end{picture}
$$
This produces the two closed oriented circle diagrams
which we already multiplied together in $H^2_2$ above.
Finally ``open'' the result to get back to a vector in $K^1_1$
by removing the left and right vertices (and the arcs
passing through them) from all the diagrams
in which the new vertices are still oriented $\down$ and
$\up$, respectively, and
mapping all other diagrams to zero:
$$
\begin{picture}(142,50)
\put(1.5,24){\oval(23,23)[t]}
\put(1.5,24){\oval(69,46)[t]}
\put(1.5,24){\oval(23,23)[b]}
\put(1.5,24){\oval(69,46)[b]}
\put(36,24){\line(-1,0){69}}
\put(-12.8,19.4){$\scriptstyle\up$}
\put(-35.8,24.1){$\scriptstyle\down$}
\put(10.2,24.1){$\scriptstyle\down$}
\put(33.2,19.4){$\scriptstyle\up$}
\put(42,21.5){$+$}
\put(91.5,24){\oval(23,23)[t]}
\put(91.5,24){\oval(69,46)[t]}
\put(91.5,24){\oval(23,23)[b]}
\put(91.5,24){\oval(69,46)[b]}
\put(126,24){\line(-1,0){69}}
\put(77.2,24.1){$\scriptstyle\down$}
\put(54.2,19.4){$\scriptstyle\up$}
\put(100.2,19.4){$\scriptstyle\up$}
\put(123.2,24.1){$\scriptstyle\down$}
\put(134,22){$\rightsquigarrow$}
\end{picture}
\begin{picture}(20,50)
\put(23.5,24){\oval(23,23)[b]}
\put(23.5,24){\oval(23,23)[t]}
\put(35,24){\line(-1,0){23}}
\put(9.2,19.4){$\scriptstyle\up$}
\put(32.2,24.1){$\scriptstyle\down$}
\put(42,22){$+\:\:0$}
\end{picture}
$$

Now we sketch the definitions of limiting versions.
First there are infinite dimensional versions
$H^\infty_m$ and $H^{\pm \infty}_m$ of $H^n_m$,
defined by
extending the number line infinitely far in both directions for $H^\infty_m$,
infinitely far to the right for $H^{+\infty}_m$, or
infinitely far to the left for $H^{-\infty}_m$.
In the diagrams parametrising the bases of these algebras,
weights have infinitely many $\up$'s and exactly $m$ $\down$'s, and
cup/cap diagrams have exactly $m$ cups/caps as before but
infinitely many rays.
Similarly, there are infinite dimensional versions
$K^{\infty}_m$, $K^{\pm \infty}_m$ and $K^\infty_\infty$
of $K^n_m$, the last of which has infinitely many $\up$'s
and infinitely many $\down$'s but finitely many cups and caps in
its diagram basis.
In all cases, the multiplications in these
infinite dimensional versions can be computed effectively
using
the generalised surgery procedure from $\S$6.
However, it is usually best to view
the algebras
$K^\infty_m, K^{\pm \infty}_m$ and $K^\infty_\infty$
as direct limits of finite dimensional $K^m_n$'s
as is explained formally in $\S$4, and then
the algebras $H^\infty_m$ and $H^{\pm \infty}_m$ emerge naturally
as subalgebras.
In fact it happens that
$H^\infty_m = K^\infty_m$ and $H^{+\infty}_m = K^{+\infty}_m$, which
is why we know that these algebras are both
quasi-hereditary and symmetric.

\vspace{2mm}
One final point.
In the main body of the article, we will adopt a slightly different notation
to this introduction, denoting the important algebras
by $H_\La$ and $K_\La$ where $\La$ is the set (or
``block'') of weights being considered.
We do this because it provides a little extra flexibility, allowing us
to incorporate additional combinatorial data into the notion of
weight
by labeling some of the vertices on the number line by $\circ$'s
and $\times$'s. These extra vertices are ``dummies'' that
play absolutely no role in
this article
but they play a key role in Parts III and IV: the extra labels encode
the combinatorics of singular blocks
on the category $\mathcal O$ and general linear supergroup sides.
There are obvious isomorphisms (``delete all $\circ$'s and $\times$'s'') between
the algebras $H_\La$ and $K_\La$ studied later on and the algebras described in
this introduction,
so we hope this causes no confusion.

\vspace{1mm}
\noindent
{\em Notation.}
In the rest of the article we work over
a fixed ground field $\C$, and
gradings mean $\Z$-gradings.

\section{Diagrams}\label{sdiagrams}

In this section we introduce the basic combinatorial
notions of weights, blocks, cup diagrams and cap diagrams, and circle diagrams.
We then explain how cup, cap and circle diagrams can
be ``oriented'' by decorating them with an additional weight.

\phantomsubsection{Weights and blocks}
By a {\em number line} we mean a horizontal line containing some {\em vertices}
indexed by a set of consecutive integers in increasing order from left to
right. A {\em weight} is a diagram obtained by labeling each of the vertices
on such a number line by $\circ$ (nought), $\times$ (cross), $\down$ (down) or
$\up$ (up) in such a way that outside of some finite subset of the vertices it
is impossible to find two (not necessarily neighbouring) vertices that are
labeled by $\down \:\up$ in that order from left to right.

A weight $\la$ is {\em bounded} if there are only finitely many
vertices. In that case, we let $\la^*$ denote the weight obtained from $\la$ by
reversing the orientations on all the vertices, i.e. changing all $\down$'s to
$\up$'s and $\up$'s to $\down$'s. In any case, we let $\la^\curvearrowleft$
denote the weight obtained from $\la$ by rotating the diagram through
$180^\circ$. For example:
$$
% scriptstyles at (1,x) x = -68,-45,-22,1,24,47,70,93,...
% \times is -0.2
% \circ is +1
% \up is +0.3 and down 2.6
% \down is +0.3 and up 2.1
\begin{picture}(30,19)
\put(-95.5,1){$\la = $}
\put(-64.8,3){\line(1,0){184}}
\put(-68.2,1){$\scriptstyle\times$}
\put(-44.8,0.2){$\circ$}
\put(46.8,1){$\scriptstyle\times$}
\put(93.2,0.2){$\circ$}
\put(-21.7,-1.6){$\scriptstyle\up$}
\put(70.3,-1.6){$\scriptstyle\up$}
\put(1.3,3.1){$\scriptstyle\down$}
\put(116.3,3.1){$\scriptstyle\down$}
\put(24.3,3.1){$\scriptstyle\down$}
\end{picture}
$$
$$
% scriptstyles at (1,x) x = -68,-45,-22,1,24,47,70,93,...
% \times is -0.2
% \circ is +1
% \up is +0.3 and down 2.6
% \down is +0.3 and up 2.1
\begin{picture}(30,19)
\put(-100,1){$\la^* = $}
\put(-64.8,3){\line(1,0){184}}
\put(-68.2,1){$\scriptstyle\times$}
\put(-44.8,0.2){$\circ$}
\put(46.8,1){$\scriptstyle\times$}
\put(93.2,0.2){$\circ$}
\put(-21.7,3.1){$\scriptstyle\down$}
\put(70.3,3.1){$\scriptstyle\down$}
\put(1.3,-1.6){$\scriptstyle\up$}
\put(116.3,-1.6){$\scriptstyle\up$}
\put(24.3,-1.6){$\scriptstyle\up$}
\end{picture}
$$
$$
% scriptstyles at (1,x) x = -68,-45,-22,1,24,47,70,93,...
% \times is -0.2
% \circ is +1
% \up is +0.3 and down 2.6
% \down is +0.3 and up 2.1
\begin{picture}(30,19)
\put(-104.3,1){$\la^\curvearrowleft = $} \put(-64.8,3){\line(1,0){184}}
\put(116.4,1){$\scriptstyle\times$}
\put(93.2,.2){$\circ$}
\put(70.3,3.1){$\scriptstyle\down$} \put(47.3,-1.6){$\scriptstyle\up$}
\put(24.3,-1.6){$\scriptstyle\up$} \put(1,1){$\scriptstyle\times$}
\put(-21.7,3.1){$\scriptstyle\down$}
\put(-44.8,0.2){$\circ$}
\put(-67.7,-1.6){$\scriptstyle\up$}
\end{picture}
$$
\vspace{1mm}

The {\em Bruhat order}
on the set of weights
is the partial order $\leq$ generated by the basic operation of
swapping a $\down$ and an $\up$; getting bigger in the Bruhat order means
$\down$'s move to the right.
We define an equivalence relation $\sim$ on the set of all
weights by declaring that $\la \sim \mu$
if $\mu$ can be obtained from $\la$ by permuting
$\down$'s and $\up$'s (and doing nothing to $\circ$'s and $\times$'s).
In other words, $\la \sim \mu$ if and only if there exists a weight
$\nu$
such that $\la \leq \nu \geq \mu$ in the Bruhat order.
We refer to the $\sim$-equivalence classes of weights as {\em blocks}.

\phantomsubsection{Cup and cap diagrams}
A {\em cup diagram} is a diagram obtained by attaching rays and finitely many
cups to some subset of the vertices on a number line, so that rays
join vertices down to infinity, cups are lower semi-circles joining pairs of
vertices, and no rays and/or cups intersect. Here is an example:
$$
% scriptstyles at (1,x) x = -68,-45,-22,1,24,47,70,93,...
% \times is -0.2
% \circ is +1
% \up is +0.3 and down 2.6
% \down is +0.3 and up 2.1
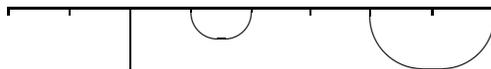
\begin{picture}(30,28)
\put(-65,23){\line(1,0){184.2}}
\put(-64.8,23){\line(0,-1){2.5}}
\put(-41.8,23){\line(0,-1){2.5}}
\put(95.5,23){\line(0,-1){2.5}}
\put(49.5,23){\line(0,-1){2.5}}
\put(95.5,23){\oval(47.2,46)[b]}
\put(15.6,23){\oval(23,23)[b]} \put(-18.8,23){\line(0,-1){23}}
\end{picture}
$$
Two cup diagrams are {\em equal} if their underlying number lines are the same
and all their rays and cups are incident with exactly the same vertices
(regardless of the precise shapes of the cups).
A vertex that is not at the end
of some cup or ray will be called a {\em free vertex},
indicated on the diagram by a vertical mark.
% and clearly displayed as $\bullet$.

The notion of a {\em cap diagram} is entirely similar, replacing cups with
caps, i.e. upper semi-circles, and drawing rays up to infinity rather than
down. For example, the following is a cap diagram:
$$
% scriptstyles at (1,x) x = -68,-45,-22,1,24,47,70,93,...
% \times is -0.2
% \circ is +1
% \up is +0.3 and down 2.6
% \down is +0.3 and up 2.1
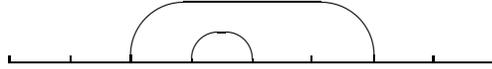
\begin{picture}(30,28)
\put(-65,1){\line(1,0){184.2}}
%\put(-66.8,-2.6)\free
%\put(-43.8,-2.6)\free
%\put(47.2,-2.6)\free
%\put(93.2,-2.6)\free
\put(95.5,1){\line(0,1){2.5}}
\put(49.5,1){\line(0,1){2.5}}
\put(-64.8,1){\line(0,1){2.5}}
\put(-41.8,1){\line(0,1){2.5}}
\put(27,1){\oval(92,46)[t]}
\put(15.5,1){\oval(23,23)[t]} \put(119,1){\line(0,1){23}}
\end{picture}
$$
A {\em closed} {cup or cap diagram} means one containing no rays. We denote the
mirror image of a cup or cap diagram $c$ in the number line by $c^*$, which is
a cap or cup diagram, respectively. Also $c^\curvearrowleft$ denotes the
diagram obtained by rotating $c$ through $180^\circ$.

If $c$ is a cup diagram and $\la$ is a weight with the same underlying
number line,
we can glue $c$ under $\la$ to obtain a new diagram
 denoted $c \la$.
We call $c \la$ an {\em oriented cup diagram} if
\begin{itemize}
\item
each free vertex in $c$
is labeled
either $\circ$ or $\times$
 in the weight $\la$;
\item
the vertices at the ends of each cup in $c$ are
labeled by exactly one $\down$ and one $\up$ in the weight $\la$;
\item
each vertex  at the top of a ray in $c$ is labeled
either $\down$ or $\up$ in the weight $\la$;
\item it is impossible to find
two rays in $c$ whose top vertices
are labeled
$\down\: \up$ in that order from left to right in the weight $\la$.
\end{itemize}
For example,
the first of
following diagrams is oriented, the other three are not (and note we now omit the vertical marks for free vertices):
$$
% scriptstyles at (1,x) x = -68,-45,-22,1,24,47,70,93,...
% \times is -0.2
% \circ is +1
% \up is +0.3 and down 2.6
% \down is +0.3 and up 2.1
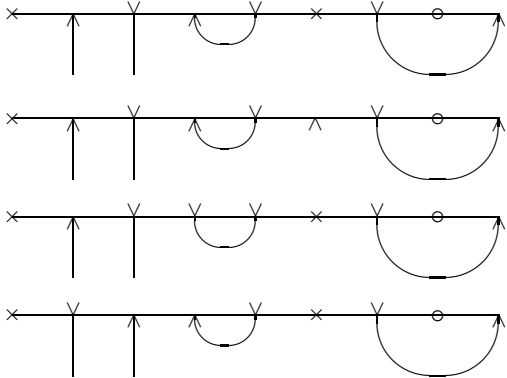
\begin{picture}(30,34)
\put(-64.8,23){\line(1,0){184}}
\put(-68.2,21){$\scriptstyle\times$}
\put(46.9,21){$\scriptstyle\times$}
\put(93.2,20.2){$\circ$}
\put(96,23){\oval(46,46)[b]}
\put(15.5,23){\oval(23,23)[b]}
\put(-19,23){\line(0,-1){23}}
\put(-42,23){\line(0,-1){23}}
\put(-44.7,18.4){$\scriptstyle\up$}
\put(-21.7,23.1){$\scriptstyle\down$}
\put(70.3,23.1){$\scriptstyle\down$}
\put(1.3,18.4){$\scriptstyle\up$}
\put(116.3,18.4){$\scriptstyle\up$}
\put(24.3,23.1){$\scriptstyle\down$}
\end{picture}
$$
$$
% scriptstyles at (1,x) x = -68,-45,-22,1,24,47,70,93,...
% \times is -0.2
% \circ is +1
% \up is +0.3 and down 2.6
% \down is +0.3 and up 2.1
\begin{picture}(30,34)
\put(-64.8,23){\line(1,0){184}}
\put(-68.2,21){$\scriptstyle\times$}
\put(46.9,18.4){$\scriptstyle\up$}
\put(93.2,20.2){$\circ$}
\put(96,23){\oval(46,46)[b]}
\put(15.5,23){\oval(23,23)[b]}
\put(-19,23){\line(0,-1){23}}
\put(-42,23){\line(0,-1){23}}
\put(-44.7,18.4){$\scriptstyle\up$}
\put(-21.7,23.1){$\scriptstyle\down$}
\put(70.3,23.1){$\scriptstyle\down$}
\put(1.3,18.4){$\scriptstyle\up$}
\put(116.3,18.4){$\scriptstyle\up$}
\put(24.3,23.1){$\scriptstyle\down$}
\end{picture}
$$
$$
% scriptstyles at (1,x) x = -68,-45,-22,1,24,47,70,93,...
% \times is -0.2
% \circ is +1
% \up is +0.3 and down 2.6
% \down is +0.3 and up 2.1
\begin{picture}(30,34)
\put(-64.8,23){\line(1,0){184}}
\put(-68.2,21){$\scriptstyle\times$}
\put(46.9,21){$\scriptstyle\times$}
\put(93.2,20.2){$\circ$}
\put(96,23){\oval(46,46)[b]}
\put(15.5,23){\oval(23,23)[b]}
\put(-19,23){\line(0,-1){23}}
\put(-42,23){\line(0,-1){23}}
\put(-44.7,18.4){$\scriptstyle\up$}
\put(-21.7,23.1){$\scriptstyle\down$}
\put(70.3,23.1){$\scriptstyle\down$}
\put(1.3,23.1){$\scriptstyle\down$}
\put(24.3,23.1){$\scriptstyle\down$}
\put(116.3,18.4){$\scriptstyle\up$}
\end{picture}
$$
$$
% scriptstyles at (1,x) x = -68,-45,-22,1,24,47,70,93,...
% \times is -0.2
% \circ is +1
% \up is +0.3 and down 2.6
% \down is +0.3 and up 2.1
\begin{picture}(30,34)
\put(-64.8,23){\line(1,0){184}}
\put(-68.2,21){$\scriptstyle\times$}
\put(46.9,21){$\scriptstyle\times$}
\put(93.2,20.2){$\circ$}
\put(96,23){\oval(46,46)[b]}
\put(15.5,23){\oval(23,23)[b]}
\put(-19,23){\line(0,-1){23}}
\put(-42,23){\line(0,-1){23}}
\put(-44.7,23.1){$\scriptstyle\down$}
\put(-21.7,18.4){$\scriptstyle\up$}
\put(70.3,23.1){$\scriptstyle\down$}
\put(1.3,18.4){$\scriptstyle\up$}
\put(116.3,18.4){$\scriptstyle\up$}
\put(24.3,23.1){$\scriptstyle\down$}
\end{picture}
$$

Instead, if $c$ is a cap diagram and $\la$ is a weight
with the same underlying number line,
we can glue $\la$ under $c$ to obtain a diagram
denoted $\la c$.
We call $\la c$ an {\em oriented cap diagram} if $c^* \la$ is an oriented
cup diagram according to the previous definition.
The following is an example of an oriented cap diagram:
$$
% scriptstyles at (1,x) x = -68,-45,-22,1,24,47,70,93,...
% \times is -0.2
% \circ is +1
% \up is +0.3 and down 2.6
% \down is +0.3 and up 2.1
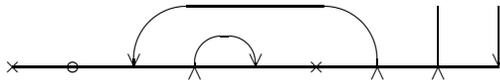
\begin{picture}(30,30)
\put(-64.8,3){\line(1,0){184}}
\put(-68.2,1){$\scriptstyle\times$}
\put(-44.8,.2){$\circ$}
\put(46.8,1){$\scriptstyle\times$}
\put(27,3){\oval(92,46)[t]}
\put(15.5,3){\oval(23,23)[t]}
\put(119,3){\line(0,1){23}}
\put(96,3){\line(0,1){23}}
\put(70.3,-1.6){$\scriptstyle\up$}
\put(-21.7,3.1){$\scriptstyle\down$}
\put(1.3,-1.6){$\scriptstyle\up$}
\put(24.3,3.1){$\scriptstyle\down$}
\put(116.3,3.1){$\scriptstyle\down$}
\put(93.3,-1.6){$\scriptstyle\up$}
\end{picture}
$$
We describe a
cup or a cap in any such oriented diagram
as {\em anti-clockwise} or {\em clockwise} according to whether its
leftmost vertex is labeled
$\down$ or $\up$ (equivalently, its rightmost vertex is labeled
$\up$ or $\down$). For example, in the above oriented cap diagram,
the small cap is clockwise and the big cap is anti-clockwise.

\phantomsubsection{Circle diagrams}
By a {\em circle diagram} we mean a diagram of the form
$ab$ obtained by gluing a cup diagram $a$ underneath a cap diagram $b$,
assuming that $a$ and $b$ have the same underlying number lines
and their free vertices are in all the same positions.
For example:
$$
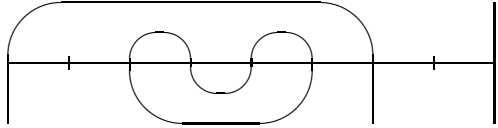
\begin{picture}(30,48)
\put(-65.2,23){\line(1,0){184.3}}
\put(-42,20.5){\line(0,1){5}}
\put(96,20.5){\line(0,1){5}}
\put(15.5,23){\oval(69,46)[b]} \put(15.5,23){\oval(23,23)[b]}
\put(38.5,23){\oval(23,23)[t]} \put(4,23){\oval(138,46)[t]}
\put(-7.5,23){\oval(23,23)[t]} \put(-65,23){\line(0,-1){23}}
\put(119,23){\line(0,-1){23}} \put(73,23){\line(0,-1){23}}
\put(119,23){\line(0,1){23}}
\end{picture}
$$
Any such diagram is a union of {\em circles}
and {\em lines} (both viewed up to homeomorphism), e.g. the above circle diagram contains one circle and two lines.
We call it a {\em closed} circle diagram
if it involves only circles and no lines.
Also a {\em small circle} means a circle consisting of just one cap and one cap.

An {\em oriented circle diagram}
means a diagram $a \la b$
obtained by gluing an oriented cup diagram $a\la$
underneath an oriented cap diagram $\la b$.
Here is an example:
$$
% scriptstyles at (1,x) x = -68,-45,-22,1,24,47,70,93,...
% \times is -0.2
% \circ is +1
% \up is +0.3 and down 2.6
% \down is +0.3 and up 2.1
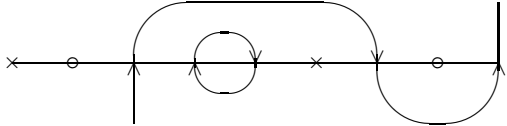
\begin{picture}(30,49)
\put(-64.8,23){\line(1,0){184}} \put(-68.2,21){$\scriptstyle\times$}
\put(-44.8,20.2){$\circ$}
\put(46.8,21){$\scriptstyle\times$}
\put(93.2,20.2){$\circ$}
\put(27,23){\oval(92,46)[t]}
\put(15.5,23){\oval(23,23)[t]} \put(119,23){\line(0,1){23}}
\put(96,23){\oval(46,46)[b]} \put(15.5,23){\oval(23,23)[b]}
\put(-19,23){\line(0,-1){23}} \put(-21.7,18.4){$\scriptstyle\up$}
\put(1.3,18.4){$\scriptstyle\up$} \put(24.3,23.1){$\scriptstyle\down$}
\put(70.3,23.1){$\scriptstyle\down$} \put(116.3,18.4){$\scriptstyle\up$}
\end{picture}
$$
We refer to a circle in an oriented circle diagram as
\begin{itemize}
\item an {\em anti-clockwise circle} or
a {circle of \em type $1$}
if the leftmost vertex on the circle is labeled by $\down$
 (equivalently, the rightmost vertex is labeled by $\up$);
\item a {\em clockwise circle} or
a {circle of \em type $x$}
if the leftmost vertex on the circle is labeled by $\up$ (equivalently, the rightmost vertex  is labeled by $\down$).
\end{itemize}

\phantomsubsection{Degrees}
The {\em degree} of an oriented cup diagram $a \la$ or
an oriented cap diagram $\la b$ means the total number of clockwise
cups and caps that it contains.
The {\em degree} of an oriented circle diagram $a \la b$ is defined by
\begin{equation}\label{dp}
\deg(a \la b) := \deg(a \la) + \deg(\la b).
\end{equation}
Equivalently, $\deg(a \la b)$ is the sum of the degrees of the
circles and lines in the diagram, where
the degree of a single circle or line in an oriented circle diagram
means
its total number of clockwise cups and caps.
For example, the oriented circle diagram displayed
at the end of the previous
subsection is of degree $3$.

\begin{Lemma}\label{deglem}
The degree of an anti-clockwise circle in an oriented circle diagram
is one less than
the total number of caps that it contains.
The degree of a clockwise circle is one more than
the total number of caps that it contains.
\end{Lemma}

\begin{proof}
Proceed by induction on the number of caps contained in the circle.
For the base case,  circles with just one cap are {small circles}
and
the lemma is clear:
$$
\begin{picture}(0,30)
\put(18.2,11.6){$\deg\Bigg($}
\put(81.2,11.6){$\Bigg) = 2.$}
\put(49.2,9.4){$\scriptstyle\up$}
\put(72.2,14.1){$\scriptstyle\down$}
\put(63.5,14){\oval(23,23)[b]}
\put(63.5,14){\oval(23,23)[t]}
\put(85,14){\line(-1,0){43}}

\put(-109.8,11.6){$\deg\Bigg($}
\put(-46.8,11.6){$\Bigg) = 0,$}
\put(-63.5,14){\oval(23,23)[b]}
\put(-63.5,14){\oval(23,23)[t]}
\put(-85,14){\line(1,0){43}}
\put(-54.8,9.4){$\scriptstyle\up$}
\put(-77.8,14.1){$\scriptstyle\down$}
\end{picture}
$$
For the induction step, remove a kink using one of the following
straightening rules:
$$
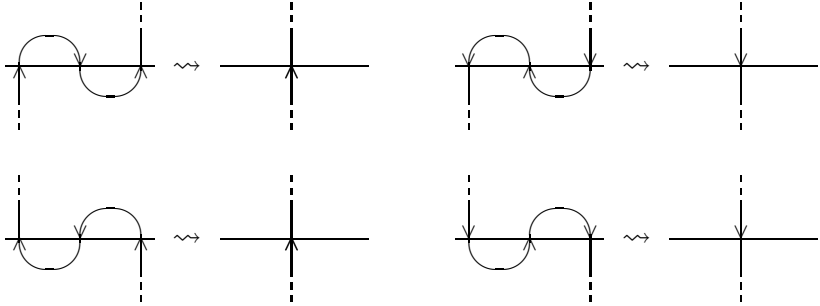
\begin{picture}(0,52)
\put(-153,21){\line(1,0){56}}
\put(-148,21){\line(0,-1){12}}
\put(-102,21){\line(0,1){12}}
\dashline{3}(-102,33)(-102,45)
\dashline{3}(-148,9)(-148,-3)
\put(-136.5,21){\oval(23,23)[t]}
\put(-113.5,21){\oval(23,23)[b]}
\put(-150.8,16.6){$\scriptstyle\up$}
\put(-127.8,21.3){$\scriptstyle\down$}
\put(-104.8,16.6){$\scriptstyle\up$}
\put(-90,18.5){$\rightsquigarrow$}
\put(-72,21){\line(1,0){56}}
\put(-45,21){\line(0,1){12}}
\put(-45,21){\line(0,-1){12}}
\put(-47.8,16.6){$\scriptstyle\up$}
\dashline{3}(-45,33)(-45,45)
\dashline{3}(-45,9)(-45,-3)

\put(17,21){\line(1,0){56}}
\put(22,21){\line(0,-1){12}}
\put(68,21){\line(0,1){12}}
\dashline{3}(68,33)(68,45)
\dashline{3}(22,9)(22,-3)
\put(33.5,21){\oval(23,23)[t]}
\put(56.5,21){\oval(23,23)[b]}
\put(19.2,21.3){$\scriptstyle\down$}
\put(42.2,16.6){$\scriptstyle\up$}
\put(65.2,21.3){$\scriptstyle\down$}
\put(80,18.5){$\rightsquigarrow$}
\put(122.2,21.3){$\scriptstyle\down$}
\put(125,21){\line(0,1){12}}
\put(125,21){\line(0,-1){12}}
\put(-47.8,16.6){$\scriptstyle\up$}
\dashline{3}(125,33)(125,45)
\dashline{3}(125,9)(125,-3)
\put(98,21){\line(1,0){56}}
\end{picture}
$$
$$
\begin{picture}(0,62)
\put(-153,21){\line(1,0){56}}
\put(-148,21){\line(0,1){12}}
\put(-102,21){\line(0,-1){12}}
\dashline{3}(-148,33)(-148,45)
\dashline{3}(-102,9)(-102,-3)
\put(-113.5,21){\oval(23,23)[t]}
\put(-136.5,21){\oval(23,23)[b]}
\put(-150.8,16.6){$\scriptstyle\up$}
\put(-127.8,21.3){$\scriptstyle\down$}
\put(-104.8,16.6){$\scriptstyle\up$}
\put(-90,18.5){$\rightsquigarrow$}
\put(-72,21){\line(1,0){56}}
\put(-45,21){\line(0,1){12}}
\put(-45,21){\line(0,-1){12}}
\put(-47.8,16.6){$\scriptstyle\up$}
\dashline{3}(-45,33)(-45,45)
\dashline{3}(-45,9)(-45,-3)

\put(17,21){\line(1,0){56}}
\put(22,21){\line(0,1){12}}
\put(68,21){\line(0,-1){12}}
\dashline{3}(22,33)(22,45)
\dashline{3}(68,9)(68,-3)
\put(56.5,21){\oval(23,23)[t]}
\put(33.5,21){\oval(23,23)[b]}
\put(19.2,21.3){$\scriptstyle\down$}
\put(42.2,16.6){$\scriptstyle\up$}
\put(65.2,21.3){$\scriptstyle\down$}
\put(80,18.5){$\rightsquigarrow$}
\put(98,21){\line(1,0){56}}
\put(122.2,21.3){$\scriptstyle\down$}
\put(125,21){\line(0,1){12}}
\put(125,21){\line(0,-1){12}}
\put(-47.8,16.6){$\scriptstyle\up$}
\dashline{3}(125,33)(125,45)
\dashline{3}(125,9)(125,-3)
\end{picture}
$$
The result is a new oriented circle drawn on a number line with two
fewer vertices than before.
Moreover, the new circle is oriented in the same way as the original circle,
and it either has one less clockwise cup and one less anti-clockwise
cap, or it has
one less anti-clockwise cup and one less clockwise
cap.
 This information suffices to complete
the proof of the induction step.
\end{proof}

\phantomsubsection{The cup diagram associated to a weight}
We will often use weights to parametrise cup diagrams
via the following important construction:
for each weight $\la$ there is a unique cup diagram
denoted $\underline{\la}$ with the property that
$\underline{\la} \la$ is an oriented cup diagram
of degree $0$, i.e. all of its cups are anti-clockwise.
To compute $\underline{\la}$ in practice, find a pair of vertices
labeled
$\down\:\up$ in order from left to right that are {neighbours}
in the sense that there
are only $\circ$'s, $\times$'s or vertices already joined by cups
in between. Join these new vertices together with a
cup then repeat the process until
there are no more such
$\down\:\up$ pairs (this will happen after
finally many repetitions). Finally draw rays down to infinity
at
all the remaining $\up$'s and $\down$'s.
For example, if $\la$ is
$$
% scriptstyles at (1,x) x = -68,-45,-22,1,24,47,70,93,...
% \times is -0.2
% \circ is +1
% \up is +0.3 and down 2.6
% \down is +0.3 and up 2.1
\begin{picture}(30,10)
\put(-64.8,3){\line(1,0){184}}
\put(-67.7,3.1){$\scriptstyle\down$}
\put(-44.7,-1.6){$\scriptstyle\up$}
\put(46.8,1){$\scriptstyle\times$}
\put(93.2,.2){$\circ$}
\put(-21.7,-1.6){$\scriptstyle\up$}
\put(70.3,-1.6){$\scriptstyle\up$}
\put(1.3,3.1){$\scriptstyle\down$}
\put(116.3,3.1){$\scriptstyle\down$}
\put(24.3,3.1){$\scriptstyle\down$}
\end{picture}
$$
then $\underline{\la}$ equals
$$
% scriptstyles at (1,x) x = -68,-45,-22,1,24,47,70,93,...
% \times is -0.2
% \circ is +1
% \up is +0.3 and down 2.6
% \down is +0.3 and up 2.1
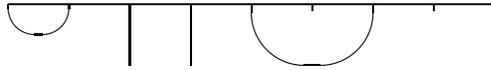
\begin{picture}(30,34)
\put(-64.8,28){\line(1,0){184}}
\put(96,28){\line(0,-1){2.5}}
\put(50,28){\line(0,-1){2.5}}
\put(-53.5,28){\oval(23,23)[b]} \put(50,28){\oval(46,46)[b]}
\put(-19,28){\line(0,-1){23}} \put(4,28){\line(0,-1){23}}
\put(119,28){\line(0,-1){23}}
\end{picture}
$$
Similarly,
$\overline{\la} := (\underline{\la})^*$
denotes the unique cap diagram such that
$\la \overline{\la}$ is an oriented cap diagram of degree $0$, i.e.
all its caps are anti-clockwise.

\begin{Lemma}\label{inj}
If $\la \sim \mu$ and $\underline{\la} = \underline{\mu}$
then $\la = \mu$.
\end{Lemma}

\begin{proof}
As $\underline{\la} = \underline{\mu}$, all the cups in the oriented cup
diagrams $\underline{\la}\la$ and $\underline{\la}\mu$ are anti-clockwise.
Hence the
labels on the vertices at the ends of each cup of $\underline{\la}$ are the
same in both $\la$ and $\mu$. On the remaining vertices, all $\up$'s must be to
the left of all $\down$'s. So since $\la \sim \mu$ all these remaining labels
are also the same in $\la$ and $\mu$.
\end{proof}

For weights $\la$ and $\mu$, we use the
notation $\mu \subset \la$ to indicate that $\mu \sim \la$ and $\underline{\mu}
\la$ is an oriented cup diagram. Sometimes we write $\la \supset \mu$ instead
of $\mu \subset \la$, so we have that $\la \supset \mu$ if $\la\sim\mu$ and
$\la \overline{\mu}$ is an oriented cap diagram. Using the first statement from
the following lemma, it is easy to check that the relation $\subset$ just
defined is reflexive and anti-symmetric; its transitive closure is the Bruhat
order $\leq$.
(Perhaps now is a good moment to recall that our cup diagrams
have only finitely many cups; this is crucial for everything
in the rest of the section.)

\begin{Lemma}\label{l2}
If $\mu \subset \la$ then $\mu \leq \la$ in the Bruhat order.
Moreover,
\begin{itemize}
\item[(i)] if $a \la$ is an oriented cup diagram then $a = \underline{\alpha}$
for a unique weight $\alpha$ with $\alpha \subset \la$;
\item[(ii)] if $\la b$ is an oriented cap diagram then $b = \overline{\beta}$
for a unique weight $\beta$ with $\la \supset \beta$;
\item[(iii)] if $a \la b$ is an oriented circle diagram
then $a = \underline{\alpha}$ and $b = \overline{\beta}$ for unique weights
$\alpha,\beta$ with $\alpha \subset \la \supset \beta$.
\end{itemize}
\end{Lemma}

\begin{proof}
Let $a \la$ be an oriented cup diagram.
We claim that $a = \underline{\alpha}$ for some weight $\alpha$ such that
$\alpha \leq \la$.
To see this,
reverse the orientations
on the vertices at the ends of each clockwise cup in the diagram
$a \la$ to obtain a new
oriented cup diagram involving only anti-clockwise cups.
This new diagram is of the form $a \alpha$ for some
weight $\alpha$, and since it is of degree $0$ we necessarily have that
$a = \underline{\alpha}$.
Whenever a clockwise cup is converted to an anti-clockwise
cup during this process, a pair of vertices labeled
$\up \:\down$ change to
$\down\:\up$, so the weight gets smaller in the Bruhat ordering.
Hence we have that $\alpha \leq \la$ and the claim follows.

If we assume now $\mu \subset \la$ and apply the claim from the previous paragraph
to $a := \underline{\mu}$,
we get that $\underline{\mu} = \underline{\alpha}$ for some weight
$\alpha$ with $\alpha \leq \la$.
Since $\mu \sim \la \sim \alpha$ we actually have that $\mu = \alpha$
by Lemma~\ref{inj}.
Hence $\mu \subset \la$ implies $\mu \leq \la$ as required.
Next we prove (i). If $a \la$ is any oriented cup diagram,
the claim in the previous paragraph shows
that $a = \underline{\alpha}$ for some $\alpha \subset \la$.
For the uniqueness, suppose also that $a = \underline{\beta}$ for $\beta
\subset
\la$. Then $\alpha \sim \la \sim \beta$ and $\underline{\alpha} =
\underline{\beta}$,
hence $\alpha = \beta$ by Lemma~\ref{inj}.
The proof of (ii) is similar. Then (iii) follows immediately from (i) and (ii).
\end{proof}

\phantomsubsection{Defect}
For a weight $\la$, the {\it defect} $\defect(\la)\in\{0,1,2,\dots\}$
means the (necessarily finite)
number of cups in the cup diagram $\underline{\la}$. The defect of a block $\La$ is then
\begin{equation}
\label{blockdefect}
\defect(\La) := \sup\{\defect(\la)\:|\:\la \in \La\}\in\{0,1,2,\dots\}\cup\{\infty\}.
\end{equation}
It is easy to see that $\defect(\La)$ is simply
the number of $\up$'s or the number of $\down$'s in the
weights belonging to $\La$, whichever is smaller.

\begin{Lemma}\label{notmany}
Let $\La$ be a block and fix $\la \in \La$.
\begin{itemize}
\item[(i)] The number of weights $\mu$ such that $\mu \supset \la$
is equal to $2^{\defect(\la)}$.
\item[(ii)] The number of weights $\mu$ such that
$\mu \subset \la$ is finite
if and only if $\defect(\La) < \infty$.
\item[(iii)]
For each fixed degree $j \in \Z$,
the number of weights
$\mu$ such that $\mu \subset \la$
and $\deg(\underline{\mu} \la) = j$ is finite
(even if $\defect(\La) = \infty$).
\end{itemize}
\end{Lemma}

\begin{proof}
(i) Start from the oriented cap diagram $\la\overline{\la}$. Choose some subset
$S$ of the caps in this diagram and reverse the orientations of the vertices at
the ends of each cap in $S$. The result is an oriented cap diagram of the form
$\mu\overline{\la}$ for some weight $\mu$ in the same block as $\la$.
Moreover, every weight $\mu$ with $\mu \supset \la$ is obtained in this way
from some choice of $S$. It remains to observe that there are
$2^{\defect(\la)}$ different choices for $S$ because $\overline{\la}$ has
$\defect(\la)$ caps.

(ii) In view of Lemma~\ref{l2}(i), the number of weights $\mu$ such that $\mu
\subset \la$ is the same as the number of cup diagrams $a$ such that $a \la$ is
an oriented cup diagram. If $\defect(\La) = \infty$ then $\la$ has infinitely
many $\up$'s and $\down$'s. So there exist oriented cup diagrams of the form $a
\la$ having an arbitrary large number of cups. Conversely, if $\defect(\La) <
\infty$ then either $\la$ has finitely many $\up$'s or finitely many $\down$'s
(or both). In the former case, suppose $a \la$ is an oriented cup diagram. At
each of the finitely many vertices labeled $\up$ in $\la$, there must be
either a ray or the leftmost end of a cup or the rightmost end of a cup in $a$.
As rays are not allowed to cross cups,
the cup diagram $a$ is uniquely determined by $\la$ and this additional
information, so there are only finitely many possibilities for $a$ (at most three to the power the number of $\up$'s in $\la$).
The argument when $\la$ has finitely many $\down$'s is analogous.

(iii) In view of (ii), we only need to consider the case that $\defect(\La) =
\infty$. We need to show that are only finitely many cup diagrams $a$ such that
$a \la$ is oriented and $\deg(a \la) = j$. To see this, note that far to the
left of $\la$ there are only $\up$'s and far to the right there are only
$\down$'s. So if $a \la$ is oriented the $\up$'s far to the left can only be at
the end of a ray or at the leftmost end of a clockwise cup, and similarly the
$\down$'s far to the right can only be at the end of a ray or at the rightmost
end of a clockwise cup. Given a fixed bound $j$ on the total number of
clockwise cups in $a$, there are still only finitely many possible $a$'s.
\end{proof}

\section{Cellularity of Khovanov's diagram algebra}\label{skhovanov}

Throughout the section, assume that $\La$ is a
{\em Khovanov block} of {rank} $n \geq 0$,
 meaning that it is a block
consisting of bounded weights
that have exactly $n$
$\up$'s and $n$ $\down$'s.
We are going to associate
a certain
algebra $H_\La$ to $\La$.
The construction of this algebra is due
to Khovanov \cite[$\S$2.4]{K2}, but we will explain it in detail
since our viewpoint is rather more combinatorial in nature than
Khovanov's. Using the new viewpoint, we then prove that $H_\La$ is
cellular.

\phantomsubsection{The underlying graded vector space}
Define $\La^{\!\circ}$ to be the subset of
$\La$ consisting of all weights
with the property that at each vertex
the number of $\down$'s
labeling vertices to the left or equal to the given vertex
is greater than or equal to
the number of $\up$'s
labeling these vertices.
In other words,
\begin{align*}
\La^{\!\circ}
&= \{\alpha \in \La\:|\:\underline{\alpha}
\text{ is a closed cup diagram}\}\\
&= \{\beta \in \La\:|\:\overline{\beta}
\text{ is a closed cap diagram}\}.
\end{align*}
We note that $|\La^{\!\circ}| = \frac{1}{n+1}\binom{2n}{n}$, the $n$th Catalan
number; see e.g. \cite[ex.6.19]{Stanley}.
Let $H_\La$ denote the vector space with basis
\begin{equation}\label{basis1}
\left\{(\underline{\alpha}\la \overline{\beta}) \:\big|\:
\text{for all }\alpha \in \La^{\!\circ}, \la \in \La,
\beta \in \La^{\!\circ}\text{ such that }
\alpha \subset \la \supset \beta\right\}.
\end{equation}
Equivalently, by Lemma~\ref{l2}(iii), the
basis of $H_\La$ is
\begin{equation}\label{equiv1}
\{(a \la b)\:|\:
\text{for all closed oriented circle diagrams $a \la b$ with $\la \in \La$}
\}.
\end{equation}
We put a grading on $H_\La$ by defining the
degree of the basis vector $( a \la b )$
to be $\deg(a \la b)$; see (\ref{dp}).
Writing $e_\la  :=
( \underline{\la}\la\overline{\la} )$ for short,
the vectors
\begin{equation}\label{ela1}
\left\{e_\la \:|\:\la \in \La^{\!\circ}\right\}
\end{equation}
form a basis for the degree $0$ component of the graded vector space $H_\La$.

\phantomsubsection{Multiplication}
Now we define a multiplication
making the graded vector space $H_\La$ into a
positively graded associative algebra.
Take two closed oriented circle diagrams
$a \la b$ and $c \mu d$
with $\la,\mu \in \La$.
If $b^* \neq c$, we simply define the product
$( a \la b )( c \mu d )$ to be zero.
Now assume that $b^* = c$.
Draw the diagram $a \la b$ underneath the diagram
$c \mu d$ to obtain a new sort of diagram
with two number lines and a symmetric middle section containing $b$ under $c$.
Then iterate the {\em surgery procedure} explained in the next paragraph
to convert this diagram into a disjoint union of diagrams none of which
have any cups or caps left in their
middle sections.
Finally, collapse each of the resulting diagrams by
identifying their top and bottom number lines
to obtain a disjoint union
of some new closed oriented circle diagrams.
The product $( a \la b ) ( c \mu d )$
is then defined to be the sum of the corresponding basis vectors of $H_\La$.

\phantomsubsection{The surgery procedure} Choose a symmetric pair of a
cup and a cap in the middle section of a diagram that can be connected without
crossings. Forgetting orientations for a while,
cut open the cup and the cap then stitch the loose ends together
to form
a pair of vertical line segments:
$$
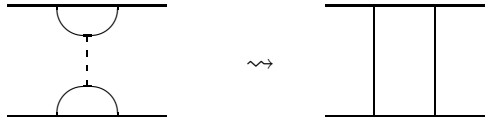
\begin{picture}(150,45)
\put(-30,-2){\line(1,0){60}}
\put(-30,40){\line(1,0){60}}
\put(0,40){\oval(23,23)[b]}
\put(0,-2){\oval(23,23)[t]}
\dashline{2.5}(0,9.5)(0,28.5)
\put(60,15){$\rightsquigarrow$}
\put(90,-2){\line(1,0){60}}
\put(90,40){\line(1,0){60}}
\put(108.5,40){\line(0,-1){42}}
\put(131.5,40){\line(0,-1){42}}
\end{picture}
$$
This surgery either combines two circles into one or splits one circle into two. Finally re-orient the new circle(s) according to the orientation of
the old circle(s) using the following rules (remembering
$1=\text{anti-clockwise}$, $x=\text{clockwise}$)
\begin{align}
1 \otimes 1 \mapsto 1,\label{mult}
\qquad
1 \otimes x \mapsto x,
\qquad
&x \otimes 1 \mapsto x,
\qquad
x \otimes x \mapsto 0,\\
1 \mapsto 1 \otimes x + x \otimes 1,
\qquad
&x \mapsto x \otimes x.\label{comult}
\end{align}
to obtain a disjoint union of zero, one or two new oriented diagrams
replacing the old diagram.
For instance, the rule $1 \otimes 1 \mapsto
1$ here indicates that two anti-clockwise circles transform to one
anti-clockwise circle. The rule $1 \mapsto 1 \otimes x + x \otimes 1$ indicates
that one anti-clockwise circle transforms to give a disjoint union of two
different diagrams each having one anti-clockwise circle and one clockwise
circle. Finally the rule $x \otimes x \mapsto 0$ indicates that two clockwise
circles produce no diagram at all at the end of the surgery procedure.

\phantomsubsection{Khovanov's definition via TQFTs}
In order to explain why the above construction makes
$H_\La$ into a well-defined associative algebra, we need to
explain in a little more detail the relationship to
Khovanov's original setup. Let
$\operatorname{2-Cob}$ be the cobordism category whose objects are closed
oriented $1$-manifolds and whose morphisms are homeomorphism classes of
oriented $2$-cobordisms. Let $\operatorname{Vec}$ be the category of vector
spaces over $\C$. The commutative Frobenius algebra $\C[x] / (x^2)$ whose
multiplication and comultiplication are defined by the formulae
(\ref{mult})--(\ref{comult}) defines a 2-dimensional TQFT, that is, a symmetric monoidal functor
\begin{equation}\label{tqft}
\mathcal F: \operatorname{2-Cob}
\rightarrow \operatorname{Vec}.
\end{equation}
We refer the reader to \cite{Kock} for the detailed construction of $\mathcal F$, just noting here
that it takes a disjoint union
$X_1 \sqcup\cdots\sqcup X_r$
of $r$ copies of $S^1$
to the tensor product of $r$ copies
of the vector space $\C[x] / (x^2)$.
The monomials $t_1 \otimes\cdots\otimes t_r$
with each $t_i \in \{1,x\}$ give a canonical
basis for the latter vector space.

Now suppose we are given a closed circle diagram
$\underline{\alpha} \overline{\beta}$ with $\alpha,\beta \in \La^{\!\circ}$.
Letting $r$ denote
the number of circles in this diagram,
we pick a homeomorphism from $\underline{\alpha}
\overline{\beta}$ to a disjoint union $X_1 \sqcup \cdots \sqcup X_r$ of $r$
copies of $S^1$.
This choice induces an isomorphism between
the vector space $\mathcal F(\underline{\alpha} \overline{\beta})$
and the tensor product of $r$ copies of $\C[x] / (x^2)$.
Now the basis $\{t_1 \otimes \cdots \otimes t_r\:|\:t_1,\dots,t_r \in \{1,x\}\}$
for the latter tensor product pulls back to give us a basis
for $\mathcal F(\underline{\alpha} \overline{\beta})$ too. It is natural to
parametrise this basis simply by decorating the $i$th circle
in the diagram $\underline{\alpha}\overline{\beta}$
by the additional {\em type} $t_i \in \{1,x\}$; equivalently, we orient the
$i$th circle
so that it is an anti-clockwise or clockwise circle
according to whether $t_i = 1$ or $t_i = x$.
In this way, the canonical basis for the vector space
$\mathcal{F}(\underline{\alpha} \overline{\beta})$ is parametrised by the
set of closed oriented circle diagrams
$$
\{\underline{\alpha}\lambda \overline{\beta}
\:|\:\text{for all $\la \in \La$ such that $\alpha \subset \la \supset \beta$}
\}.
$$
Now it makes sense to identify the vector space $H_\La$ from \eqref{basis1}
with the vector space
\begin{equation}\label{kid}
\bigoplus_{\alpha,\beta \in \La^{\!\circ}}
\mathcal F(\underline{\alpha}\overline{\beta}),
\end{equation}
since both spaces have distinguished bases parameterised by the same set.

Under this identification,
the multiplication on $H_\La$ described above
corresponds
to the direct sum of linear maps
$$
m_{\alpha,\beta,\gamma,\delta}:\quad \mathcal
F(\underline{\alpha}\overline{\beta}) \otimes \mathcal
F(\underline{\gamma}\overline{\delta}) \rightarrow \mathcal
F(\underline{\alpha}\overline{\delta})
$$
for all $\alpha,\beta,\gamma,\delta \in \La^{\!\circ}$ defined as follows. If
$\beta \neq \gamma$, then $m_{\alpha,\beta,\gamma,\delta} = 0$. If $\beta =
\gamma$, then $m_{\alpha,\beta,\gamma,\delta}$ is the linear map arising by
applying the functor $\mathcal F$ to the cobordism
$\underline{\alpha}\overline{\beta} \:\sqcup\:
\underline{\gamma}\overline{\delta} \rightarrow
\underline{\alpha}\overline{\delta}$ constructed by composing ``pair of pants''
shaped cobordisms going from two to one or from one to two circles lifting the sequence of surgery procedures explained above.
The advantage of this more topological
point of view is that it is now clear  that
the multiplication
is well defined independent of the order that
the surgery procedures are performed, since the resulting
cobordism is independent of the choice of order up to homeomorphism;
see \cite[p.679]{K2} for a careful justification of the latter statement.
In a similar fashion, the proof that the multiplication now defined
is associative
reduces to checking that two corresponding
cobordisms are homeomorphic, something
which is apparent from their construction.

We still need to prove the multiplication preserves the grading, a
fact which is also proved on \cite[p.680]{K2} from the topological point of
view. We have to verify that each surgery procedure
preserves
the total number of clockwise cups and caps in a diagram.
In view of Lemma~\ref{deglem}, this amounts to showing that
$$
\#(\text{caps}) - \#(\text{anti-clockwise circles}) + \#(\text{clockwise circles})
$$
is the same at the end of each surgery procedure
as it was at the start.
When two circles go to one circle,
there is one less cap, one less anti-clockwise circle and the same number of clockwise circles
at the end compared to the beginning.
When one circle goes to two circles,
there is one less cap, the same number of anti-clockwise
circles and one more clockwise circle
at the end compared to the beginning.
The desired invariance follows directly.

\phantomsubsection{Cellular algebra structure}
We have now constructed the graded associative algebra $H_\La$.
It is clear from the definition that
\begin{eqnarray*}
e_\al ( a \la b )=
\begin{cases}
( a \la b )&\text{if $\underline{\al}=a$,}\\
0&\text{otherwise,}
\end{cases}
&\quad&
( a \la b ) e_\be
=
\begin{cases}
( a \la b ) &\text{if $b = \overline{\be}$,}\\
0&\text{otherwise,}
\end{cases}
\end{eqnarray*}
for $\al,\be \in \La^{\!\circ}$ and any basis vector $(a \la b)\in H_\La$.
This implies that the vectors
$\{e_\al\:|\:\al \in \La^{\!\circ}\}$ are mutually orthogonal
idempotents whose sum is the identity in  $H_\La$.
The resulting direct sum decomposition
\begin{equation}
H_\La = \bigoplus_{\alpha,\beta \in \La^{\!\circ}} e_\alpha\; H_\La\; e_\beta
\end{equation}
is just the decomposition from (\ref{kid}) in different notation, and the
summand $e_\alpha H_\La e_\beta$ has basis $\left\{(\underline{\alpha} \la
\overline{\beta})\:\big|\:\text{for all } \la \in \La\text{ such that } \alpha
\subset \la \supset \beta\right\}$.
Observe also that the linear maps $*:H_\La
\rightarrow H_\La$ and $\curvearrowleft:H_\La \rightarrow
H_{\La^\curvearrowleft}$ defined  by
\begin{align}
\label{antiaut}
( a \la b )^* &:= ( b^* \la a^* ),\\
( a \la b
)^\curvearrowleft &:= ( b^\curvearrowleft \la^\curvearrowleft a^\curvearrowleft
)
\end{align}
are graded algebra anti-isomorphisms, where
$\La^\curvearrowleft:=\{\la^\curvearrowleft\mid\la\in\La\}$.

\begin{Theorem}\label{cell}
Let $( a \la b )$ and
$( c \mu d )$ be basis vectors of $H_\La$.
Then,
$$
( a \la b )
( c \mu d ) =
\left\{
\begin{array}{ll}
0&\text{if $b \neq c^*$,}\\
s_{a \la b}(\mu) ( a \mu d ) + (\dagger)&\text{if $b = c^*$
and $a \mu$ is oriented,}\\
(\dagger)&\text{otherwise,}
\end{array}
\right.
$$
where
\begin{itemize}
\item[(i)]
$(\dagger)$ denotes a linear combination of basis vectors
of $H_\La$ of the form $( a \nu d)$
for $\nu > \mu$;
\item[(ii)]
the scalar $s_{a \la b}(\mu) \in \{0,1\}$ depends only on $a \la b$ and $\mu$,
but not on $d$;
\item[(iii)] if $b = \overline{\la}=c^*$
and
$a\mu$ is oriented
then $s_{a \la b}(\mu) = 1$.
\end{itemize}
\end{Theorem}

\begin{proof}
We already know that $( a \la b ) ( c \mu d ) = 0$ if $b \neq c^*$, so assume
$b=c^*$ from now on. Consider a single iteration of the surgery procedure in the
algorithm computing $( a \la b ) ( c \mu d )$. Let $\tau$ be the {top
weight} (the weight on the top number line) of this diagram
at the start of the
surgery procedure. We claim that
\begin{itemize}
\item[(1)]
the top weight of each diagram obtained at the end of
the surgery procedure is greater than or equal to $\tau$ in the Bruhat order;
\item[(2)]
the total number of
diagrams produced with top weight equal to $\tau$
is either zero or one, independent of the cap diagram $d$;
\item[(3)]
if $a \tau$ is oriented
and the cap to be cut is
anti-clockwise then exactly one
diagram is produced with top weight equal to $\tau$.
\end{itemize}
Before proving these things, let us explain how the theorem follows.
By applying (1) repeatedly, starting with $\tau = \mu$ at the first step,
it follows that
$(a \la b) (c \mu d)$ is a linear combination of
$(a \nu d)$'s for $\nu \geq \mu$.
Assuming $a \mu$ is oriented,
% (so that $(a \mu d)$ is a well defined
%basis vector of $H_\La$),
(2) applied repeatedly
implies that the coefficient
$s_{a \la b}(\mu)$ of the basis vector
$(a \mu d)$ in the product is zero or one
independent of the cap diagram $d$.
This proves (i) and (ii).
If in addition
$b = \overline{\la}$ (so that {\em all} caps to be cut
are anti-clockwise) then
(3) applied repeatedly implies that  $s_{a \la b}(\mu)=1$,
as in (iii).
Therefore, it just remains to prove (1)--(3).
To do this, we analyse three different situations
depending on the orientations of the cap and the cup to be cut
(both of which clearly do not depend on $d$).

\noindent
{\em Case one: the cap to be cut is clockwise.}
In this case, we show that the top weight of each diagram
obtained at the end of the surgery procedure is strictly greater than $\tau$,
from which (1)--(3) follow immediately.
If the circle containing the cap to be
cut is anti-clockwise, the cap is necessarily a concave part of the
circumference of that circle, and we are in one of three basic situations
represented by the following pictures:
$$
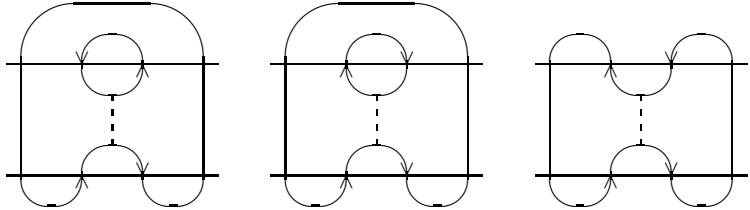
\begin{picture}(0,79)
\put(-140,10){\line(1,0){80}}
\put(-140,52){\line(1,0){80}}
\put(-100,52){\oval(23,23)[b]}
\put(-100,52){\oval(23,23)[t]}

\put(-100,52){\oval(69,46)[t]}
\put(-100,10){\oval(23,23)[t]}
\put(-77,10){\oval(23,23)[b]}
\put(-123,10){\oval(23,23)[b]}
\dashline{2.5}(-100,21.5)(-100,40.5)
\put(-134.5,10){\line(0,1){42}}
\put(-65.5,10){\line(0,1){42}}
\put(-114.3,5.4){$\scriptstyle\up$}
\put(-91.3,10.1){$\scriptstyle\down$}
\put(-91.3,47.4){$\scriptstyle\up$}
\put(-114.3,52.1){$\scriptstyle\down$}
\end{picture}
\begin{picture}(0,77)
\put(-40,10){\line(1,0){80}}
\put(-40,52){\line(1,0){80}}
\put(0,52){\oval(23,23)[b]}
\put(0,52){\oval(23,23)[t]}
\put(0,52){\oval(69,46)[t]}
\put(0,10){\oval(23,23)[t]}
\put(23,10){\oval(23,23)[b]}
\put(-23,10){\oval(23,23)[b]}
\dashline{2.5}(0,21.5)(0,40.5)
\put(-34.5,10){\line(0,1){42}}
\put(34.5,10){\line(0,1){42}}
\put(-14.3,5.4){$\scriptstyle\up$}
\put(8.7,10.1){$\scriptstyle\down$}
\put(-14.3,47.4){$\scriptstyle\up$}
\put(8.7,52.1){$\scriptstyle\down$}
\end{picture}
\begin{picture}(0,77)
\put(60,10){\line(1,0){80}}
\put(60,52){\line(1,0){80}}
\put(100,52){\oval(23,23)[b]}
\put(77,52){\oval(23,23)[t]}
\put(123,52){\oval(23,23)[t]}
\put(100,10){\oval(23,23)[t]}
\put(123,10){\oval(23,23)[b]}
\put(77,10){\oval(23,23)[b]}
\dashline{2.5}(100,21.5)(100,40.5)
\put(65.5,10){\line(0,1){42}}
\put(134.5,10){\line(0,1){42}}
\put(85.7,5.4){$\scriptstyle\up$}
\put(108.7,10.1){$\scriptstyle\down$}
\put(85.7,47.4){$\scriptstyle\up$}
\put(108.7,52.1){$\scriptstyle\down$}
\end{picture}
$$
(We stress that these and subsequent pictures
should be interpreted only up to homeomorphism;
in particular the circles represented in the pictures may well
cross both number lines many more
times than indicated.)
In the first of these cases, the orientation of each
vertex lying on the anti-clockwise circle containing the cup to be cut
gets reversed when re-orienting,
hence the top weight gets strictly larger in the Bruhat order
as required.
For the second case, the orientation of each vertex
from the anti-clockwise circle containing the cap to be cut
gets reversed. Since this circle necessarily crosses the top number line, we see again that the top weight increases.
For the final case, two anti-clockwise circles
are produced when the initial surgery is performed.
Then these circles are re-oriented so that orientation of all
vertices lying on one or other of these
anti-clockwise circles gets switched,
so the top weight again goes up in the Bruhat order.
Instead, assume the circle containing the cap to be cut is clockwise,
so the cap is a convex part of the circumference.
There are then four basic situations:
$$
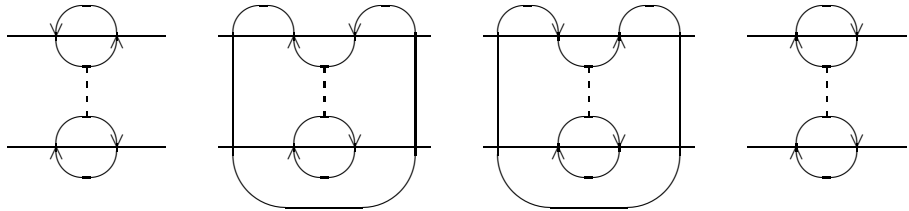
\begin{picture}(0,77)
\put(-170,20){\line(1,0){60}}
\put(-170,62){\line(1,0){60}}
\put(-140,20){\oval(23,23)[t]}
\put(-140,62){\oval(23,23)[t]}
\put(-140,20){\oval(23,23)[b]}
\put(-140,62){\oval(23,23)[b]}
\dashline{2.5}(-140,31.5)(-140,50.5)
\put(-154.3,15.4){$\scriptstyle\up$}
\put(-131.3,20.1){$\scriptstyle\down$}
\put(-131.3,57.4){$\scriptstyle\up$}
\put(-154.3,62.1){$\scriptstyle\down$}
\end{picture}
\begin{picture}(0,77)
\put(-90,20){\line(1,0){80}}
\put(-90,62){\line(1,0){80}}
\put(-50,20){\oval(23,23)[t]}
\put(-50,20){\oval(23,23)[b]}
\put(-50,20){\oval(69,46)[b]}
\put(-50,62){\oval(23,23)[b]}
\put(-27,62){\oval(23,23)[t]}
\put(-73,62){\oval(23,23)[t]}
\dashline{2.5}(-50,31.5)(-50,50.5)
\put(-84.5,20){\line(0,1){42}}
\put(-15.5,20){\line(0,1){42}}
\put(-64.3,15.4){$\scriptstyle\up$}
\put(-41.3,20.1){$\scriptstyle\down$}
\put(-64.3,57.4){$\scriptstyle\up$}
\put(-41.3,62.1){$\scriptstyle\down$}
\end{picture}
\begin{picture}(0,77)
\put(10,20){\line(1,0){80}}
\put(10,62){\line(1,0){80}}
\put(50,20){\oval(23,23)[t]}
\put(50,20){\oval(23,23)[b]}
\put(50,20){\oval(69,46)[b]}
\put(50,62){\oval(23,23)[b]}
\put(73,62){\oval(23,23)[t]}
\put(27,62){\oval(23,23)[t]}
\dashline{2.5}(50,31.5)(50,50.5)
\put(15.5,20){\line(0,1){42}}
\put(84.5,20){\line(0,1){42}}
\put(35.7,15.4){$\scriptstyle\up$}
\put(58.7,20.1){$\scriptstyle\down$}
\put(58.7,57.4){$\scriptstyle\up$}
\put(35.7,62.1){$\scriptstyle\down$}
\end{picture}
\begin{picture}(0,77)
\put(110,20){\line(1,0){60}}
\put(110,62){\line(1,0){60}}
\put(140,20){\oval(23,23)[t]}
\put(140,62){\oval(23,23)[t]}
\put(140,20){\oval(23,23)[b]}
\put(140,62){\oval(23,23)[b]}
\dashline{2.5}(140,31.5)(140,50.5)
\put(125.7,15.4){$\scriptstyle\up$}
\put(148.7,20.1){$\scriptstyle\down$}
\put(125.7,57.4){$\scriptstyle\up$}
\put(148.7,62.1){$\scriptstyle\down$}
\end{picture}
$$
In the first two of these, the orientation of
every vertex lying on the anti-clockwise circle containing
the cup to be cut
gets switched, so the top weight gets bigger.
The last two involve the rule $x \otimes x \mapsto 0$,
so there is nothing to check.

\noindent
{\em Case two: both the cap and the cup to be cut are anti-clockwise.}
Here we show there is exactly one new
diagram produced having top weight equal to $\tau$, together possibly with one
other diagram with a strictly larger top weight. This is again enough to
establish (1)--(3). There are five basic configurations to consider:
$$
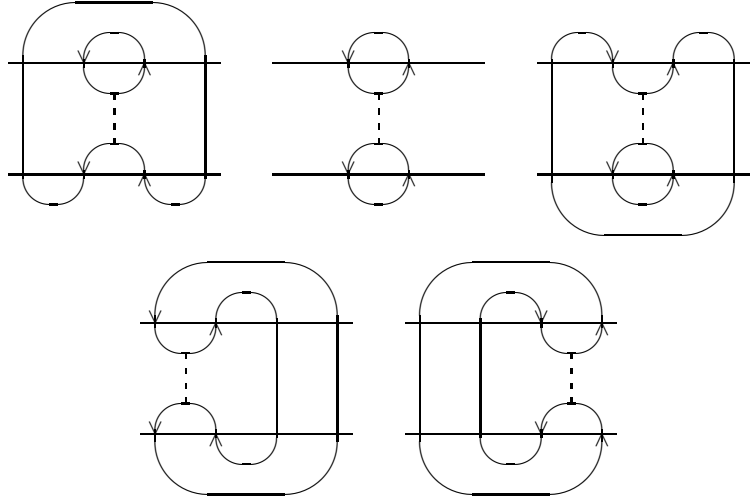
\begin{picture}(0,95)
\put(-40,22){\line(1,0){80}}
\put(-40,64){\line(1,0){80}}
\put(0,64){\oval(23,23)[t]}
\put(0,22){\oval(23,23)[t]}
\put(0,64){\oval(23,23)[b]}
\put(0,22){\oval(23,23)[b]}
\dashline{2.5}(0,33.5)(0,52.5)
\put(8.7,17.4){$\scriptstyle\up$}
\put(-14.3,22.1){$\scriptstyle\down$}
\put(8.7,59.4){$\scriptstyle\up$}
\put(-14.3,64.1){$\scriptstyle\down$}
\end{picture}
\begin{picture}(0,95)
\put(60,22){\line(1,0){80}}
\put(60,64){\line(1,0){80}}
 \put(100,22){\oval(23,23)[b]}
\put(100,22){\oval(23,23)[t]}
\put(100,22){\oval(69,46)[b]}
\put(100,64){\oval(23,23)[b]}
\put(123,64){\oval(23,23)[t]}
\put(77,64){\oval(23,23)[t]}
\dashline{2.5}(100,33.5)(100,52.5)
\put(65.5,22){\line(0,1){42}}
\put(134.5,22){\line(0,1){42}}
\put(108.7,17.4){$\scriptstyle\up$}
\put(85.7,22.1){$\scriptstyle\down$}
\put(108.7,59.4){$\scriptstyle\up$}
\put(85.7,64.1){$\scriptstyle\down$}
\end{picture}
\begin{picture}(0,95)
\put(-140,22){\line(1,0){80}}
\put(-140,64){\line(1,0){80}}
\put(-100,64){\oval(23,23)[b]}
\put(-100,64){\oval(23,23)[t]}
\put(-100,64){\oval(69,46)[t]}
\put(-100,22){\oval(23,23)[t]}
\put(-77,22){\oval(23,23)[b]}
\put(-123,22){\oval(23,23)[b]}
\dashline{2.5}(-100,33.5)(-100,52.5)
\put(-134.5,22){\line(0,1){42}}
\put(-65.5,22){\line(0,1){42}}
\put(-91.3,17.4){$\scriptstyle\up$}
\put(-114.3,22.1){$\scriptstyle\down$}
\put(-91.3,59.4){$\scriptstyle\up$}
\put(-114.3,64.1){$\scriptstyle\down$}
\end{picture}
$$
$$
\begin{picture}(0,95)
\put(10,22){\line(1,0){80}}
\put(10,64){\line(1,0){80}}
\put(73,64){\oval(23,23)[b]}
\put(50,64){\oval(69,46)[t]}
\put(50,64){\oval(23,23)[t]}
\put(50,22){\oval(69,46)[b]}
\put(73,22){\oval(23,23)[t]}
\put(50,22){\oval(23,23)[b]}
\dashline{2.5}(73,33.5)(73,52.5)
\put(15.5,22){\line(0,1){42}}
\put(38.5,22){\line(0,1){42}}
\put(81.7,17.4){$\scriptstyle\up$}
\put(58.7,22.1){$\scriptstyle\down$}
\put(81.7,59.4){$\scriptstyle\up$}
\put(58.7,64.1){$\scriptstyle\down$}
\end{picture}
\begin{picture}(0,95)
\put(-90,22){\line(1,0){80}}
\put(-90,64){\line(1,0){80}}
\put(-73,64){\oval(23,23)[b]}
\put(-50,64){\oval(69,46)[t]}
\put(-50,64){\oval(23,23)[t]}
\put(-50,22){\oval(69,46)[b]}
\put(-73,22){\oval(23,23)[t]}
\put(-50,22){\oval(23,23)[b]}
\dashline{2.5}(-73,33.5)(-73,52.5)
\put(-38.5,22){\line(0,1){42}}
\put(-15.5,22){\line(0,1){42}}
\put(-64.3,17.4){$\scriptstyle\up$}
\put(-87.3,22.1){$\scriptstyle\down$}
\put(-64.3,59.4){$\scriptstyle\up$}
\put(-87.3,64.1){$\scriptstyle\down$}
\end{picture}
$$
In the first three of these configurations,
no vertices change orientation during the procedure,
so we get precisely one diagram at the end, and it still has
top weight equal to $\tau$.
For both
the remaining two, we obtain a pair of nested circles,
 one anti-clockwise and the other clockwise, when the surgery is performed.
This produces one new diagram with the top weight equal to
$\tau$.
Re-orienting the diagram produces
another new
diagram obtained from this by reversing the orientations of
both of these circles. Since the outside circle is anti-clockwise
and the inside one is clockwise, the top weight of the resulting
diagram
is indeed strictly greater than $\tau$ in the Bruhat order.

\noindent
{\it Case three: the cap to be cut is anti-clockwise but the cup to be cut is clockwise.}
In this case, the outcome
depends on whether the circle containing the cap to be cut crosses the top
number line (something which is again independent of $d$) or not.
We show if it does
not cross the top number line then just one diagram, with top weight $\tau$, is produced by the surgery
procedure, and
if it does cross the top number line then we only get diagrams
with a strictly larger top weight.
This is enough to prove (1)--(2); then we will argue further in the next paragraph to
get (3).
To see these things, there are three remaining situations to be considered:
$$
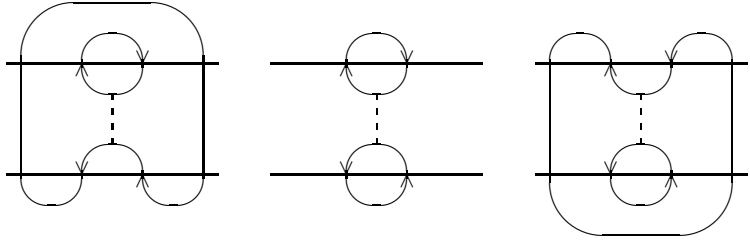
\begin{picture}(0,95)
\put(60,22){\line(1,0){80}}
\put(60,64){\line(1,0){80}}
\put(100,22){\oval(23,23)[b]}
\put(100,22){\oval(23,23)[t]}
\put(100,22){\oval(69,46)[b]}
\put(100,64){\oval(23,23)[b]}
\put(123,64){\oval(23,23)[t]}
\put(77,64){\oval(23,23)[t]}
\dashline{2.5}(100,33.5)(100,52.5)
\put(134.5,22){\line(0,1){42}}
\put(65.5,22){\line(0,1){42}}
\put(108.7,17.4){$\scriptstyle\up$}
\put(85.7,22.1){$\scriptstyle\down$}
\put(85.7,59.4){$\scriptstyle\up$}
\put(108.7,64.1){$\scriptstyle\down$}
\end{picture}
\begin{picture}(0,95)
\put(-40,22){\line(1,0){80}}
\put(-40,64){\line(1,0){80}}
\put(0,64){\oval(23,23)[t]}
\put(0,22){\oval(23,23)[t]}
\put(0,64){\oval(23,23)[b]}
\put(0,22){\oval(23,23)[b]}
\dashline{2.5}(0,33.5)(0,52.5)
\put(8.7,17.4){$\scriptstyle\up$}
\put(-14.3,22.1){$\scriptstyle\down$}
\put(-14.3,59.4){$\scriptstyle\up$}
\put(8.7,64.1){$\scriptstyle\down$}
\end{picture}
\begin{picture}(0,95)
\put(-140,22){\line(1,0){80}}
\put(-140,64){\line(1,0){80}}
\put(-100,64){\oval(23,23)[b]}
\put(-100,64){\oval(23,23)[t]}
\put(-100,64){\oval(69,46)[t]}
\put(-100,22){\oval(23,23)[t]}
\put(-77,22){\oval(23,23)[b]}
\put(-123,22){\oval(23,23)[b]}
\dashline{2.5}(-100,33.5)(-100,52.5)
\put(-134.5,22){\line(0,1){42}}
\put(-65.5,22){\line(0,1){42}}
\put(-91.3,17.4){$\scriptstyle\up$}
\put(-114.3,22.1){$\scriptstyle\down$}
\put(-114.3,59.4){$\scriptstyle\up$}
\put(-91.3,64.1){$\scriptstyle\down$}
\end{picture}
$$
The first of these (in which the bottom circle definitely crosses the top number line) produces no diagrams at all as it is a situation
when $x \otimes x \mapsto 0$.
For the second two, the orientation of every vertex from the
anti-clockwise circle containing the cap to be cut
gets switched, so we get
a diagram with strictly bigger top weight if the bottom circle
crosses the top number line, and we get a diagram with top weight equal
to $\tau$ if it does not.

We still need to complete the proof of
(3) in case three. To do that, we must show that
if $a \tau$ is oriented then the circle containing the cap to be cut
definitely does not cross the top number line, hence
by the previous paragraph exactly one diagram
is definitely produced with top weight $\tau$.
Let $\sigma$ be the bottom weight of the diagram at the start of the surgery
procedure. Also let $e$ be
the cap diagram in the middle section of the diagram, so that the middle section of the diagram looks like $e$ glued under $e^*$
with the ends of matching rays joined together to form vertical line segments.
Because the entire diagram is oriented at the start of the surgery procedure,
both the diagrams $a \sigma e$ and $e^* \tau d$ are oriented.
So under our assumption that $a \tau$ is oriented, we get that
$a \sigma e$ and $a \tau e$ are oriented circle diagrams.
Let $i$ be the integer indexing
the leftmost vertex of the cap to be cut.
This vertex has different orientations in $a\sigma e$ and $a \tau e$.
Hence all the other vertices that are
joined to the vertex $i$
have different orientations in $a \sigma e$ and $a \tau e$ too.
Finally observe that the vertices in $\sigma$ and $\tau$ at the
ends of each vertical line segment in the middle section
necessarily have the same orientation.
Hence there cannot be any such line segment in the circle
containing the cap to be cut, i.e. this circle
does not cross the top number line.
\end{proof}

\begin{Corollary}\label{ideal}
The product
$( a \la b )
( c \mu d )$ of two basis vectors of $H_\La$
is a linear combination of vectors of the form
$( a \nu d)$ for $\nu \in \La$ with $\la \leq \nu \geq \mu$.
\end{Corollary}

\begin{proof}
By Theorem~\ref{cell}(i),
$( a \la b ) ( c \mu d )$ is a linear combination
of $( a \nu d)$'s for various $\nu \geq \mu$
and
$( d^* \mu c^* ) ( b^* \la a^* )$
is a linear combination of $( d^* \nu a^*)$'s for various
$\nu \geq \la$.
Applying the anti-automorphism $*$ from (\ref{antiaut})
to the latter statement gives
that
$( a \la b ) ( c \mu d )$ is a linear combination
of $( a \nu d)$'s for various $\nu \geq \la$ too.
\end{proof}

\begin{Corollary}\label{iscell}
Khovanov's algebra $H_\La$ is a cellular algebra in the sense of Graham and
Lehrer \cite{GL} with cell datum $(\La, M, C, *)$ where
\begin{itemize}
\item[(i)]
$M(\la)$ denotes $\left\{\alpha \in \La^{\!\circ}\:|\:\alpha
\subset \la\right\}$
for each $\la \in \La$;
\item[(ii)]
$C$ is defined by setting
$C^\la_{\alpha,\beta}:=( \underline{\alpha} \la \overline{\beta} )$
for $\la \in \La$ and $\alpha,\beta \in M(\la)$;
\item[(iii)] $*$ is the anti-automorphism from (\ref{antiaut}).
\end{itemize}
\end{Corollary}

Before we prove the corollary let us first recall the relevant definitions from
\cite{GL}. A {\em cellular algebra} means an associative unital algebra $H$
together with a {\em cell datum} $(\La, M, C, *)$ such that
\begin{itemize}
\item[(1)]
$\La$ is a partially ordered set and $M(\la)$
is a finite set
for each $\la \in \La$;
\item[(2)]
$C:\dot\bigcup_{\la \in \La} M(\la)
\times M(\la) \rightarrow H,
(\alpha,\beta) \mapsto C^\la_{\alpha,\beta}$ is an injective map
whose image is a basis for $H$;
\item[(3)]
the map $*:H \rightarrow H$
is an algebra anti-automorphism such that
$(C^\la_{\alpha,\beta})^* = C_{\beta,\alpha}^\la$
for all $\la \in \La$ and $\alpha, \beta \in M(\la)$;
\item[(4)]
if $\mu \in \La$ and $\gamma, \delta \in M(\la)$
then for any $x \in H$ we have that
$$
x C_{\gamma,\delta}^\mu \equiv \sum_{\gamma' \in M(\mu)} r_x(\gamma',\gamma) C_{\gamma',\delta}^\mu
\pmod{H(> \mu)}$$
where the scalar $r_x(\gamma',\gamma)$ is independent of $\delta$ and
$H(> \mu)$ denotes the subspace of $H$
generated by $\{C_{\gamma'',\delta''}^\nu\:|\:\nu > \mu,
\gamma'',\delta'' \in M(\nu)\}$.
\end{itemize}

\begin{proof}[Proof of Corollary~\ref{iscell}]
Condition (1) is clear as $|\La^{\!\circ}|$ itself is a finite set.
Condition (2) is a
consequence of the definition (\ref{basis1}). We have already asserted that $*$
is an anti-automorphism giving (3). Finally to verify (4) it suffices to
consider the case that $x = C^\la_{\alpha,\beta}$ for some $\la \in \La$ and
$\alpha, \beta \in M(\la)$. If $\beta = \gamma$ and $\alpha \subset \mu$, then
Theorem~\ref{cell}(i)--(ii) shows that
$$
C_{\alpha,\beta}^\la C_{\gamma,\delta}^\mu \equiv
s_{\underline{\alpha} \la \overline{\beta}}(\mu) C_{\alpha, \delta}^\mu
\pmod{H_\La(> \mu)}
$$
where $s_{\underline{\alpha} \la \overline{\beta}}(\mu)$
is independent of $\delta$;
otherwise, we have that
$$
C_{\alpha, \beta}^\la C_{\gamma, \delta}^\mu \equiv 0
\pmod{H_\La(> \mu)}.
$$
Taking
$$
r_x(\gamma', \gamma)
:=
\left\{
\begin{array}{ll}
s_{\underline{\alpha}\la\overline{\beta}}(\mu)&
\text{if $\gamma' = \alpha, \beta = \gamma$ and $\alpha \subset \mu$,}\\
0&\text{otherwise,}
\end{array}\right.
$$
we deduce that (4) holds.
\end{proof}

\begin{Remark}\rm
In fact $H_\La$ is an example of a {\em graded} cellular algebra;
the degree of $\alpha \in M(\lambda)$ is the
number of clockwise cups in the diagram $\underline{\alpha} \lambda$.
The same remark applies to Corollary~\ref{iscell2} and Theorem~\ref{iscell3} below.
\end{Remark}

\begin{Remark}\rm
Assume for this remark that $n > 0$.
Khovanov showed in \cite[Proposition 32]{K2}
that the algebra $H_\La$ is a symmetric algebra; see also Theorem~\ref{salg}
below.
Since $H_\La$ is obviously not  a semisimple algebra,
it follows easily that $H_\La$ has infinite global dimension.
In particular, this implies that $H_\La$ is
definitely {\em not} a quasi-hereditary algebra
in the sense of Cline, Parshall and Scott \cite{CPS1}.
\end{Remark}

\section{The algebra $K_\La$: a larger cellular algebra}\label{sstroppel}
Fix a block $\La$ throughout the section. We are going to construct a graded
algebra $K_\La$ which, if $\La$ is a Khovanov block, we will show in Part II
is a
quasi-hereditary cover of Khovanov's algebra $H_\La$.
We split the definition of the multiplication on
$K_\La$ into two cases: the bounded case
(when $K_\La$ is finite dimensional) and the general case (when it may be
infinite dimensional).

\phantomsubsection{The underlying graded vector space} We begin by defining $K_\La$ to be the vector space with basis
\begin{equation}\label{basis2}
\left\{
( \underline{\alpha} \la \overline{\beta})\:\big|\:
\text{for all }\alpha,\la,\beta \in \La \text{ such that }
\alpha \subset \la \supset \beta\right\}.
\end{equation}
Equivalently, by Lemma~\ref{l2}(iii), this basis is the set
\begin{equation}\label{equiv2}
\left\{( a\la b)\:\big|\:
\text{for all oriented circle diagrams $a \la b$ with $\la \in \La$}
\right\}.
\end{equation}
We stress the
differences between this and Khovanov's algebra $H_\La$:
first, we now allow more general (perhaps even infinite) blocks $\La$;
second, the
basis is parametrised by arbitrary oriented circle diagrams whereas
Khovanov's algebra only involves closed ones.
We put a grading on $K_\La$
by defining the degree of the basis vector $( a \la b )$ to be $\deg(a \la b)$ as in \eqref{dp}.
For $\la \in \La$, we write simply $e_\la$ for $( \underline{\la}\la
\overline{\la} )$. The vectors
\begin{equation}\label{ela2}
\left\{e_\la \:|\:\la \in \La\right\}
\end{equation}
give a basis
for the degree $0$ component of $K_\La$.

\phantomsubsection{Multiplication in the bounded case}
We restrict our attention for a while to the situations when $\La$
consists of bounded weights.
In that case, our strategy to define the multiplication
on $K_\La$ is to first
construct an isomorphism of graded vector spaces
$\cl:K_\La \rightarrow H_\De / I_\La$ for some Khovanov block $\De$
and some homogeneous ideal $I_\La$ of $H_\De$.
Using this isomorphism, we will then simply
lift the multiplication on the quotient algebra $H_\De /
I_\La$ to the vector space $K_\La$. We already want to point out that actual computations needed in practise can be done in a quicker way;
see $\S$\ref{section6} for details.

Let $p$ be the number of
$\up$'s and $q$ be the number of $\down$'s
labeling the vertices of each weight in $\La$.
For $\la \in \La$, we define its {\em closure}
$\cl(\la)$
to be the weight obtained from $\la$ by adding
$p$ new vertices labeled by $\down$ to the left  end of the number line
and $q$ new vertices labeled by $\up$ to the right  end.
For example, the closure of
$$
% scriptstyles at (1,x) x = -68,-45,-22,1,24,47,70,93,...
% \times is -0.2
% \circ is +1
% \up is +0.3 and down 2.6
% \down is +0.3 and up 2.1
\begin{picture}(30,10)
\put(-64.8,3){\line(1,0){184}}
\put(-68.2,1){$\scriptstyle\times$}
\put(-44.8,.2){$\circ$}
\put(46.8,1){$\scriptstyle\times$}
\put(93.2,.2){$\circ$}
\put(1.3,-1.6){$\scriptstyle\up$}
\put(24.3,3.1){$\scriptstyle\down$}
\put(-21.7,-1.6){$\scriptstyle\up$}
\put(70.3,3.1){$\scriptstyle\down$}
\put(116.3,-1.6){$\scriptstyle\up$}
\end{picture}
$$
is
$$
% scriptstyles at (1,x) x = -68,-45,-22,1,24,47,70,93,...
% \times is -0.2
% \circ is +1
% \up is +0.3 and down 2.6
% \down is +0.3 and up 2.1
\begin{picture}(30,10)
\put(-133.8,3){\line(1,0){299}}
%\put(-64.8,3){\line(1,0){184}}
%\dashline{2.5}(-64.8,3)(-133.8,3)
%\dashline{2.5}(119.2,3)(165.2,3)
\put(-68.2,1){$\scriptstyle\times$}
\put(-44.8,.2){$\circ$}
\put(46.8,1){$\scriptstyle\times$}
\put(93.2,.2){$\circ$}
\put(1.3,-1.6){$\scriptstyle\up$}
\put(24.3,3.1){$\scriptstyle\down$}
\put(-21.7,-1.6){$\scriptstyle\up$}
\put(70.3,3.1){$\scriptstyle\down$}
\put(116.3,-1.6){$\scriptstyle\up$}
\put(-90.7,3.1){$\scriptstyle\down$}
\put(-113.7,3.1){$\scriptstyle\down$}
\put(-136.7,3.1){$\scriptstyle\down$}
\put(139.3,-1.6){$\scriptstyle\up$}
\put(162.3,-1.6){$\scriptstyle\up$}
\end{picture}
$$
Let $\De$ be the unique block
containing the weights $\cl(\la)$ for all $\la \in \La$.
Note that $\De$ is a Khovanov block,
so it makes sense to consider Khovanov's algebra
$H_{\De}$ from $\S$\ref{skhovanov}.
The map
\begin{equation}\label{cl1}
\cl:\La \hookrightarrow \De
\end{equation}
just defined is injective with image $\cl(\La)$ contained in
$\De^{\!\circ}$.

\begin{Lemma}\label{upperset}
If $\la \in \De \setminus \cl(\La)$
and $\mu \geq \la$ then $\mu \in \De \setminus \cl(\La)$.
\end{Lemma}

\begin{proof}
The set $\De \setminus \cl(\La)$
consists of all weights in $\De$
either with some $\up$
labeling one of its leftmost $p$ vertices or some
$\down$ labeling one of its rightmost $q$ vertices.
The lemma is obvious from this description and the definition of the
Bruhat order.
\end{proof}

Now suppose that $a$ is a cup diagram with
$a = \underline{\alpha}$ for $\alpha \in \La$.
We define its closure $\cl(a)$ by setting
\begin{equation}\label{cl2}
\cl(a) := \underline{\cl(\alpha)}.
\end{equation}
To compute $\cl(a)$ in practice, let $r$ denote the number of cups in $a$.
First add $p$ new vertices to the left end and $q$ new
vertices to the right end of the number line in $a$. Then
connect the leftmost $r$ of these new vertices to the rightmost $r$ of the new
vertices by adding $r$ nested
cups. Finally convert the leftmost $(p-r)$ rays of $a$
into cups incident with the remaining $(p-r)$ new vertices on the left hand
side and convert the rightmost $(q-r)$ rays of $a$ into cups incident with the
remaining $(q-r)$ new vertices on the right hand side. This explicit
description implies for any $\la \in \La$ that $a \la$ is an oriented cup
diagram if and only if $\cl(a \la) := \cl(a) \cl(\la)$ is a closed oriented cup
diagram. Moreover, since only anti-clockwise cups are added when passing from
$a \la$ to $\cl(a\la)$, we have that $\deg(\cl(a\la)) = \deg(a \la)$.
Similarly, if $b$ is a cap diagram with $b = \overline{\beta}$ for $\beta \in
\La$, we set
\begin{equation}\label{cl3}
\cl(b) := \overline{\cl(\beta)},
\end{equation}
so that $\cl(b)^* = \cl(b^*)$. We then get for any oriented circle diagram $a
\la b$ with $\la \in \La$ that
\begin{equation}\label{notion}
\cl(a \la b) := \cl(a) \cl(\la) \cl(b)
\end{equation}
is a closed oriented circle diagram of the same degree as $a \la b$. For
example, if $a \la b$ equals
$$
% scriptstyles at (1,x) x = -68,-45,-22,1,24,47,70,93,...
% \times is -0.2
% \circ is +1
% \up is +0.3 and down 2.6
% \down is +0.3 and up 2.1
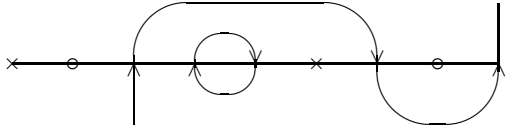
\begin{picture}(30,50)
\put(-64.8,25){\line(1,0){184}}
\put(-68.2,23){$\scriptstyle\times$}
\put(-44.8,22.2){$\circ$}
\put(119,25){\line(0,1){23}}
\put(-19,25){\line(0,-1){23}}
\put(46.8,23){$\scriptstyle\times$}
\put(93.2,22.2){$\circ$}
\put(27,25){\oval(92,46)[t]}
\put(15.5,25){\oval(23,23)[t]}
\put(96,25){\oval(46,46)[b]}
\put(15.5,25){\oval(23,23)[b]}
\put(1.3,20.4){$\scriptstyle\up$}
\put(24.3,25.1){$\scriptstyle\down$}
\put(-21.7,20.4){$\scriptstyle\up$}
\put(70.3,25.1){$\scriptstyle\down$}
\put(116.3,20.4){$\scriptstyle\up$}
\end{picture}
$$
then $\cl(a \la b)$ is the closed oriented circle diagram
$$
% scriptstyles at (1,x) x = -68,-45,-22,1,24,47,70,93,...
% \times is -0.2
% \circ is +1
% \up is +0.3 and down 2.6
% \down is +0.3 and up 2.1
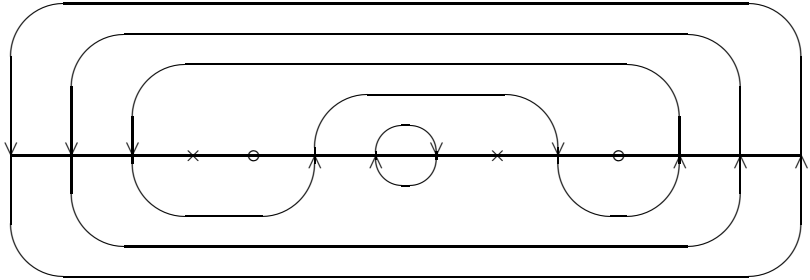
\begin{picture}(30,105)
\put(-133.8,45){\line(1,0){299}}
%\put(-64.8,45){\line(1,0){184}}
%\dashline{2.5}(-64.8,45)(-133.8,45)
%\dashline{2.5}(119.2,45)(165.2,45)
\put(-68.2,43){$\scriptstyle\times$}
\put(-44.8,42.2){$\circ$}
\put(46.8,43){$\scriptstyle\times$}
\put(93.2,42.2){$\circ$}
\put(27,45){\oval(92,46)[t]}
\put(15.5,45){\oval(23,23)[t]}
\put(96,45){\oval(46,46)[b]}
\put(15.5,45){\oval(23,23)[b]}
\put(1.3,40.4){$\scriptstyle\up$}
\put(24.3,45.1){$\scriptstyle\down$}
\put(-21.7,40.4){$\scriptstyle\up$}
\put(70.3,45.1){$\scriptstyle\down$}
\put(116.3,40.4){$\scriptstyle\up$}
\put(-90.7,45.1){$\scriptstyle\down$}
\put(-113.7,45.1){$\scriptstyle\down$}
\put(-136.7,45.1){$\scriptstyle\down$}
\put(139.3,40.4){$\scriptstyle\up$}
\put(162.3,40.4){$\scriptstyle\up$}
\put(15.5,45){\oval(207,69)[t]}
\put(-53.5,45){\oval(69,46)[b]}
\put(15.5,45){\oval(253,69)[b]}
\put(15.5,45){\oval(253,92)[t]}
\put(15.5,45){\oval(299,92)[b]}
\put(15.5,45){\oval(299,115)[t]}
\end{picture}
$$
Conversely, given a
closed oriented circle diagram with underlying weight belonging
to $\cl(\La)$, there is an obvious way to ``open'' the diagram, deleting
the leftmost $p$ and rightmost $q$ vertices in the process,
to obtain an oriented circle diagram
$a \la b$ with $\la \in \La$ such that $\cl(a \la b)$ is the initial
closed oriented circle diagram.
This proves the following lemma.

\begin{Lemma}\label{opening}
The map $a \la b \mapsto \cl(a \la b)$ is a degree preserving bijection
between the set of oriented circle diagrams with underlying weight
belonging to $\La$ and the set of closed oriented circle diagrams with
underlying weight belonging to $\cl(\La)$.
\end{Lemma}

Now we can complete the definition of the graded algebra $K_\La$
when $\La$ consists of bounded weights.
Lemma~\ref{upperset} and
Corollary~\ref{ideal} imply that the subspace
$I_\La$ of $H_{\De}$
spanned by the vectors
$$
\left\{( a \la b )\:\big|\:
\text{for all closed oriented circle diagrams
$a \la b$ with $\la \in \De \setminus \cl(\La)$}
\right\}
$$
is a two-sided ideal of $H_\De$.
In view of Lemma~\ref{opening}, the vectors
$$
\left\{( \cl(a \la b)) + I_\La\:\big|\:
\text{for all oriented circle diagrams
$a \la b$ with $\la \in \La$}\right\}
$$
give a basis for the
quotient algebra $H_\De / I_\La$.
We deduce that the map
\begin{equation}\label{closure}
\cl:K_\La \rightarrow H_\De / I_\La,\qquad
( a \la b ) \mapsto ( \cl (a \la b) ) +
I_\La
\end{equation}
is an isomorphism of graded vector spaces.
Using this, we transport the multiplication on the
algebra $H_\De / I_\La$ to the vector space $K_\La$,
to make $K_\La$ into a well-defined positively graded algebra in its own right.
In other words, to compute the product
$( a \la b ) ( c \mu d )$
of two basis vectors of $K_\La$, we first close both diagrams,
then compute their product in Khovanov's algebra $H_\De$ modulo
the ideal $I_\La$,
then apply the inverse map
$\cl^{-1}$ to get back to an element of $K_\La$.

\phantomsubsection{Multiplication in the general case}
It just remains to define the multiplication on $K_\La$ in the unbounded cases
too, by taking a direct limit of the algebras already defined. To set this up,
we need one more piece of notation. For weights $\la$ and $\mu$ we write $\la
\prec \mu$ if
\begin{itemize}
\item $\la$ is a bounded weight;
\item the set $I$
indexing the vertices of $\la$ is a subset of the set $J$ indexing
the vertices of $\mu$;
\item
the labels on each vertex of $\mu$ indexed by $I$ are the same as the labels on
the corresponding vertices of $\la$;
\item
amongst all vertices of $\mu$ indexed by the set $J \setminus I$
it is impossible to find two that are labeled
$\down\:\up$ in that order from left to right.
\end{itemize}
For example the following weights $\la$ and $\mu$
(whose vertices are indexed by $\{1,2,3\}$ and by
$\Z$, respectively)
satisfy
$\la \prec \mu$:
$$
% scriptstyles at (1,x) x = -68,-45,-22,1,24,47,70,93,...
%
\begin{picture}(30,12)
\put(4.2,3){\line(1,0){46}}
\put(24.2,.2){$\circ$}
\put(1.3,-1.4){$\scriptstyle\up$}
\put(47.3,3.1){$\scriptstyle\down$}
\end{picture}
$$
$$
% scriptstyles at (1,x) x = -68,-45,-22,1,24,47,70,93,...
%
\begin{picture}(30,12)
\put(-118,2.5){$\cdots$}
\put(-98,5){\line(1,0){223}}
%\put(-41.9,2){\line(0,1){6}}
\put(24.2,2.2){$\circ$}
\put(1.3,0.6){$\scriptstyle\up$}
\put(47.3,5.1){$\scriptstyle\down$}
\put(-21.7,0.6){$\scriptstyle\up$}
\put(-44.7,3){$\scriptstyle\times$}
\put(-67.7,0.6){$\scriptstyle\up$}
\put(-90.7,0.6){$\scriptstyle\up$}
\put(70.3,0.6){$\scriptstyle\up$}
\put(93.3,5.1){$\scriptstyle\down$}
\put(116.3,5.1){$\scriptstyle\down$}
\put(133,2.5){$\cdots$}
\end{picture}
$$
Given two blocks $\Ga$ and $\La$
we write $\Ga \prec \La$
if there exist weights $\la \in \Ga$ and $\mu \in \La$
such that $\la \prec \mu$.
In that case, for every $\la \in \Ga$ there is a unique
element denoted $\ex_\Ga^\La(\la) \in \La$,
the {\em extension} of $\la$ from $\Ga$ to $\La$,
with the property that $\la \prec \ex_\Ga^\La(\la)$.
This defines an injection
\begin{equation}\label{ex1}
\ex_\Ga^\La:\Ga \hookrightarrow \La.
\end{equation}
If $\Ga \prec \Om \prec \La$
then $\Ga \prec \La$ and
$\ex_\Ga^\La =
\ex_{\Om}^\La \circ \ex_{\Ga}^{\Om}$.
Moreover,
$\La$ is the union of the finite subsets $\ex_{\Ga}^\La(\Ga)$
for all $\Ga \prec \La$.
In other words, the maps $\ex_{\Ga}^\La$ induce a canonical
isomorphism
$\varinjlim \Ga \stackrel{\sim}{\rightarrow} \La$
between
the direct limit
of all $\Ga \prec \La$
and the set $\La$.

Now suppose that $a$ is a cup diagram with
$a = \underline{\alpha}$
and $b$ is a cap diagram with $b = \overline{\beta}$
for $\alpha,\beta \in \Ga$.
We define their extensions
$\ex_\Ga^\La(a)$ and $\ex_\Ga^\La(b)$ by setting
\begin{equation}\label{ex2}
\ex_\Ga^\La(a) := \underline{\ex_\Ga^\La(\alpha)},
\qquad
\ex_\Ga^\La(b) := \overline{\ex_\Ga^\La(\beta)}.
\end{equation}
Then for an oriented circle diagram $a \la b$
with $\la \in \Ga$, we
define its extension
to be the oriented circle diagram
\begin{equation}
\ex_\Ga^\La (a\la b)
:= \ex_\Ga^\La(a) \ex_\Ga^\La(\la) \ex_\Ga^\La(b).
\end{equation}
Informally, $\ex_\Ga^\La(a\la b)$ is obtained from the
diagram $a \la b$ by first extending the weight $\la$
to $\ex_\Ga^\La(\la)$, then extending the cup and cap diagrams
$a$ and $b$ by converting as many as possible of their
rays into anti-clockwise cups and
caps through the new vertices,
then adding new rays through any
remaining vertices labeled $\up$ or $\down$.
For instance if $\Ga$ and $\La$ are the blocks generated by the weights
$\la$ and $\mu$ from the previous example
and $a \la b$ equals
$$
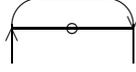
\begin{picture}(30,18)
\put(4.2,3){\line(1,0){46}}
\put(24.2,.2){$\circ$}
\put(47.3,3.1){$\scriptstyle\down$}
\put(1.3,-1.4){$\scriptstyle\up$}
\put(27.2,3){\oval(46,23)[t]}
\put(4.2,3){\line(0,-1){13}}
\put(50.2,3){\line(0,-1){13}}
\end{picture}
$$
then $\ex_\Ga^\La(a \la b)$ is
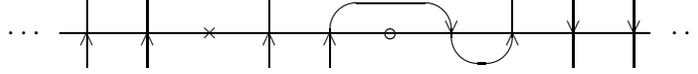
$$
% scriptstyles at (1,x) x = -68,-45,-22,1,24,47,70,93,...
%
\begin{picture}(30,26)
\put(-118,7.5){$\cdots$}
\put(-98,10){\line(1,0){223}}
%\put(-41.9,2){\line(0,1){6}}
\put(24.2,7.2){$\circ$}
\put(1.3,5.6){$\scriptstyle\up$}
\put(47.3,10.1){$\scriptstyle\down$}
\put(-21.7,5.6){$\scriptstyle\up$}
\put(-44.7,8){$\scriptstyle\times$}
\put(-67.7,5.6){$\scriptstyle\up$}
\put(-90.7,5.6){$\scriptstyle\up$}
\put(70.3,5.6){$\scriptstyle\up$}
\put(93.3,10.1){$\scriptstyle\down$}
\put(116.3,10.1){$\scriptstyle\down$}
\put(133,7.5){$\cdots$}
\put(27.2,10){\oval(46,23)[t]}
\put(61.7,10){\oval(23,23)[b]}
\put(73.2,10){\line(0,1){13}}
\put(96.2,10){\line(0,1){13}}
\put(119.2,10){\line(0,1){13}}
\put(-18.8,10){\line(0,1){13}}
\put(-64.8,10){\line(0,1){13}}
\put(-87.8,10){\line(0,1){13}}
\put(4.2,10){\line(0,-1){13}}
\put(96.2,10){\line(0,-1){13}}
\put(119.2,10){\line(0,-1){13}}
\put(-18.8,10){\line(0,-1){13}}
\put(-64.8,10){\line(0,-1){13}}
\put(-87.8,10){\line(0,-1){13}}
\end{picture}
$$
In this way, we obtain a degree preserving
injective linear map
\begin{equation}\label{extension0}
\ex_{\Ga}^{\La}:K_\Ga \hookrightarrow K_\La, \qquad ( a \la b ) \mapsto (
\ex_{\Ga}^{\La}(a\la b) ).
\end{equation}
If $\Ga \prec \Om \prec \La$,
then
$\ex_{\Ga}^{\La} =
\ex_{\Om}^{\La} \circ \ex_{\Ga}^{\Om}$,
and also $K_\La$ is the sum of the subspaces
$\ex_{\Ga}^{\La} (K_\Ga)$ for all
$\Ga \prec \La$.
So the maps $\ex_\Ga^{\La}$ induce a canonical isomorphism
of graded vector spaces
$\varinjlim K_\Ga \stackrel{\sim}{\rightarrow} K_\La$,
where the direct limit is again taken over all $\Ga \prec \La$.

\begin{Lemma}\label{ad}
Given blocks $\Ga \prec \Om \prec \La$,
the map $\ex_{\Ga}^{\Om}: K_\Ga \rightarrow K_{\Om}$
is an algebra homomorphism.
\end{Lemma}

\begin{proof}
It suffices to consider the four
cases in which $\ex_{\Ga}^{\Om}:\Ga
\rightarrow \Om$ adds one new vertex labeled either $\down$ or $\up$
to the weights in $\Ga$, either at the left or the right end.
In all four cases, for any oriented circle
diagram $a \la b$ with $\la \in \Ga$, the closure of $\ex_{\Ga}^{\Om}(a \la b)$
is obtained from the closure of $a \la b$ by adding one more
anti-clockwise circle.
More precisely, if $\ex_\Ga^\Om$ adds one vertex labeled
$\up$ to the left or one vertex labeled $\down$ to the right then $\cl(\ex_\Ga^\Om(a \la b))$ is
$\cl(a \la b)$ with one additional small anti-clockwise circle inside,
and if $\ex_\Ga^\Om$ adds one vertex labeled
$\down$ to the left or one vertex labeled $\up$ to the right then
$\cl(\ex_\Ga^\Om(a \la b))$ is
$\cl(a \la b)$ with one additional large anti-clockwise circle outside.
Since the additional circle
always sits in the same position independent of the diagram $a \la b$,
it
follows that it makes no difference to the multiplication rule, because $1
\otimes 1 \mapsto 1$ in the definition of Khovanov's algebra.
\end{proof}

In view of Lemma~\ref{ad}, the direct limit $\varinjlim K_\Ga$ is actually a
graded associative (but in general not unital) algebra. Now transport the
multiplication on $\varinjlim K_\Ga$ to the vector space $K_\La$ using the
canonical isomorphism constructed just before Lemma~\ref{ad}, to make $K_\La$
into a well-defined graded associative algebra. In other words, to compute the
product $( a \la b ) ( c \mu d )$ of two basis vectors of $K_\La$, first choose
some $\Ga \prec \La$ so that $( a \la b)$ and $( c \mu d)$ are in the image of
the map $\ex_\Ga^\La$, then ``truncate'' by applying $(\ex_\Ga^{\La})^{-1}$ and
compute the product in the finite dimensional algebra  $K_\Ga$, then finally
apply $\ex_\Ga^\La$ to extend back to an element of $K_\La$.

\phantomsubsection{Cellular algebra structure}
We have now defined the multiplication making
the graded vector space
$K_\La$ into a graded associative algebra in all cases.
The definition gives quite easily that
\begin{eqnarray*}
e_\al ( a \la b )=
\begin{cases}
( a \la b )&\text{if $\underline{\al}=a$,}\\
0&\text{otherwise,}
\end{cases}
&\quad&
( a \la b ) e_\be
=
\begin{cases}
( a \la b ) &\text{if $b = \overline{\be}$,}\\
0&\text{otherwise,}
\end{cases}
\end{eqnarray*}
for $\al,\be \in \La$ and any basis vector $( a \la b ) \in K_\La$.
So the vectors
$\{e_\al\:|\:\al \in \La\}$ give a set
of mutually orthogonal
idempotents in $K_\La$.
If $|\La| < \infty$, the sum of these idempotents is the identity element
of the finite dimensional algebra $K_\La$.
In general, $K_\La$ might not be a unital algebra but we still have at least that
\begin{equation}\label{still}
K_\La = \bigoplus_{\alpha,\beta \in \La}
e_\alpha K_\La e_\beta.
\end{equation}
The summand $e_\alpha K_\La e_\beta$
of $K_\La$ has basis $$
\left\{( \underline{\alpha} \la
\overline{\beta})\:\big|\:\text{for all } \lambda \in \La \text{ such that
}\alpha\subset\la \supset \beta\right\}.$$
This set is finite thanks to
Lemma~\ref{notmany}(i), so each $e_\alpha K_\La e_\beta$ is finite dimensional
even if $K_\La$ itself is infinite dimensional.
Note also for $\Ga \prec \La$ that the
map $\ex_\Ga^\La: K_\Ga \rightarrow K_\La$
sends the idempotent $e_\ga \in K_\Ga$ for $\ga \in \Ga$
to the idempotent $e_{\la} \in K_\La$
where $\la := \ex_\Ga^\La(\ga)$.
Observe finally
that there are graded algebra anti-isomorphisms
\begin{align}
\label{star}
*:K_\La &\rightarrow K_\La,
\qquad\qquad
(a \la b) \mapsto (b^* \la a^*),\\
\curvearrowleft:K_\La &\rightarrow K_{\La^{\curvearrowleft}},
\qquad\;\quad
(a \la b)
\mapsto (b^\curvearrowleft \la^\curvearrowleft a^\curvearrowleft),
\label{curve}
\end{align}
just like \eqref{antiaut}.
%is a graded algebra anti-automorphism, as is the map $\curvearrowleft:K_\La
%\rightarrow K_{\La^\curvearrowleft}$ defined by
%\begin{equation}\label{hash}(a \la b)^\curvearrowleft := (b^\curvearrowleft
%\la^\curvearrowleft a^\curvearrowleft),
%\end{equation}
%where
%\begin{equation}
%\La^\curvearrowleft := \{\la^\curvearrowleft\:|\:\la \in \La\}.
%\end{equation}
Composing $*$ and $\curvearrowleft$, we deduce that $K_\La \cong
K_{\La^\curvearrowleft}$ as graded algebras.

\begin{Theorem}\label{cell2}
Let $( a \la b )$ and
$( c \mu d )$ be basis vectors of $K_\La$.
Then,
$$
( a \la b )
( c \mu d ) =
\left\{
\begin{array}{ll}
0&\text{if $b \neq c^*$,}\\
s_{a \la b}(\mu) ( a \mu d ) + (\dagger)&\text{if $b = c^*$
and $a \mu$ is oriented,}\\
(\dagger)&\text{otherwise,}
\end{array}
\right.
$$
where
\begin{itemize}
\item[(i)]
$(\dagger)$ denotes a linear combination of basis vectors
of $K_\La$ of the form $( a \nu d)$
for $\nu > \mu$;
\item[(ii)]
the scalar $s_{a \la b}(\mu) \in \{0,1\}$ depends only on
$a \la b$ and $\mu$ (but not on $d$);
\item[(iii)] if $b = \overline{\la}=c^*$
and
$a\mu$ is oriented
then $s_{a \la b}(\mu) = 1$.
\end{itemize}
\end{Theorem}

\begin{proof}
Assume to start with
that the weights in $\La$ are bounded.
We have already seen that
$( a \la b ) ( c \mu d ) = 0$ if
$b^* \neq c$, so assume that $b^* =c$.
In that case, $\cl(b)^* = \cl(b^*) = \cl(c)$,
so we get by
Theorem~\ref{cell}(i) that
$( \cl(a\la b) )
( \cl(c\mu d) )
=
s_{\cl(a\la b)}({\cl(\mu)})( \cl(a\mu d) ) + (\dagger)$
if $\cl(a\mu)$ is oriented, and
$( \cl(a\la b) )
( \cl(c\mu d) )
=
(\dagger)$
otherwise,
where $(\dagger)$
denotes a linear combination of
$(\cl(a) \nu' \cl(d))$'s for $\nu' > \cl(\mu)$.
If $\nu' > \cl(\mu)$ then we either
have that $\nu' = \cl(\nu)$ for some $\nu > \mu$
or $(\cl(a) \nu' \cl(d)) \in I_\La$.
Moreover,
$\cl(a\mu)$ is oriented if and only if $a \mu$ is oriented.
So applying the map $(\cl)^{-1}$, we have shown that
$$
( a\la b )
( c\mu d )
=
\left\{
\begin{array}{ll}
s_{a \la b}({\mu})( a\mu d ) + (\dagger)&\text{if $a \mu$ is oriented,}\\
(\dagger)&\text{otherwise,}
\end{array}\right.
$$
where $s_{a \la b}(\mu) := s_{\cl(a \la b)}(\cl(\mu)) \in \{0,1\}$
and $(\dagger)$ denotes a linear combination of
basis vectors of $K_\La$ of the form $(a \nu d)$ for $\nu > \mu$.

Since $s_{a \la b}(\mu)$ depends only on $a \la b$ and $\mu$
by Theorem~\ref{cell}(ii),
this is what we wanted for (i)--(ii).
Moreover,
if $b = \overline{\la} = c^*$ and $a \mu$ is oriented as in (iii),
then we get that $s_{a \la b}(\mu) = 1$
by Theorem~\ref{cell}(iii).

It just remains to
pass from the bounded case to the unbounded case by taking direct limits.
This is done by a similar argument to the preceeding one and causes no additional problems.
\end{proof}

\begin{Corollary}\label{ideal2}
The product
$( a \la b )
( c \mu d )$ of two basis vectors of $K_\La$
is a linear combination of vectors of the form
$( a \nu d)$ for $\nu \in \La$ with $\la \leq \nu \geq \mu$.
\end{Corollary}

\begin{proof}
Argue as in the proof of Corollary~\ref{ideal}.
\end{proof}

In the next corollary, we need to exclude the situations when
$|\La| = \infty$, but only because
cellular algebras are unital algebras by definition.

\begin{Corollary}\label{iscell2}
Assuming $|\La| < \infty$,
the algebra $K_\La$
is a cellular algebra in the sense of Graham and Lehrer with cell datum
$(\La, M, C, *)$ where
\begin{itemize}
\item[(i)]
$M(\la)$ denotes $\left\{\alpha \in \La\:|\:\alpha\subset \la\right\}$
for each $\la \in \La$;
\item[(ii)]
$C$ is defined by setting
$C^\la_{\alpha,\beta}:=( \underline{\alpha} \la \overline{\beta} )$
for $\la \in \La$ and $\alpha,\beta \in M(\la)$;
\item[(iii)]
$*$ is the anti-automorphism from (\ref{star}).
\end{itemize}
\end{Corollary}

\begin{proof}
Argue as in the proof of Corollary~\ref{iscell}.
\end{proof}

\section{The quasi-hereditary structure of $K_\La$}\label{sdecomposition}

Now we begin the study of the representation theory of the algebras
$K_\La$ for any block $\La$ in a systematic way.
We define cell modules for $K_\La$, giving
explicit closed formulae for the $q$-decomposition numbers
which simultaneously describe  composition multiplicities of
cell modules and  cell filtration
multiplicities of projective indecomposable modules.
We deduce in almost all cases that the category $\rep{K_\La}$
of finite dimensional graded $K_\La$-modules is a graded highest weight
category.

\phantomsubsection{Graded modules}
We begin by fixing some conventions for working with graded modules
over a locally unital graded
algebra $K= \bigoplus_{i \in \Z} K_i$.
Here
{\em locally unital} means
that $K$ is an associative (but not necessarily unital) algebra
with a given system
$\{e_\la\:|\:\la\in \La\}$ of mutually orthogonal idempotents
such that
$$
K = \bigoplus_{\al,\beta \in \La} e_\al  K e_\be.
$$
Of course if $K$ is finite dimensional, all this is just saying
that $K$ is a unital algebra with identity element
$\sum_{\la \in \La} e_\la$.
(Alternatively, a locally unital algebra such as
$K$ could be viewed as a small category with objects
equal to
$\La$ and $\hom(\la,\mu) := e_\mu K e_\la$.)
By a {\em $K$-module} we always mean a left $K$-module
$M$ such that
$$
M = \bigoplus_{\la \in \La} e_\la M.
$$
If $K$ is finite dimensional, this is just saying that $M$ is a
unital module in the usual sense.

If $M=\bigoplus_{j \in \Z} M_j$ is a {\em graded} $K$-module,
i.e. $K_i M_j \subseteq M_{i+j}$,
then we
write $M\langle j \rangle$ for the same module but with new grading defined by
$M\langle j \rangle_i := M_{i-j}$. For graded modules $M$ and $N$, we define
\begin{equation}\label{homdef}
\hom_{K}(M, N) := \bigoplus_{j \in \Z} \hom_{K}(M,N)_j
\end{equation}
where $\hom_{K}(M,N)_j$ denotes all homogeneous $K$-module
homomorphisms of degree $j$, meaning that they map $M_i$ into $N_{i+j}$ for
each $i \in \Z$.
Let $\mod{K}$ be the category of all graded $K$-modules $M =
\bigoplus_{j \in \Z} M_j$ that are locally finite dimensional and bounded
below, i.e. each $M_j$ is finite dimensional and $M_j = 0$ for $j \ll 0$.
Morphisms in $\mod{K}$ mean homogeneous $K$-module homomorphisms of
degree $0$, so that $\mod{K}$ is an abelian category. The space $\hom_{K}(M,N)_j$
can be also described as all morphisms from $M\langle j\rangle$ to $N$ in the
category $\mod{K}$.

\phantomsubsection{\boldmath Irreducible and projective $K_\La$-modules}
Now fix an arbitrary block $\La$
and note by (\ref{still}) that $K_\La$ is locally unital with respect to the system
of idempotents $\{e_\la\:|\:\la \in \La\}$.
Let $(K_\La)_{> 0}$ be the sum of all
components of the graded algebra $K_\La$ of strictly positive degree, so
\begin{equation}\label{degzero}
K_\La / (K_\La)_{> 0} \cong \bigoplus_{\la \in \La} \C
\end{equation}
as an algebra, with a basis given by the images of
all the idempotents $e_\la$. The image of $e_\la$
spans a one dimensional graded $K_\La$-module which we denote by $L(\la)$.
Thus, $L(\la)$ is a copy of the field concentrated in degree $0$, and $(a \mu
b)\in K_\La$ acts on $L(\la)$ as multiplication by $1$ if $a \mu b =
\underline{\la} \la \overline{\la}$, or as zero otherwise. The modules
\begin{equation}\label{simples}
\{L(\la)\langle j \rangle\:|\:\la \in \La, j\in \Z\}
\end{equation}
give a complete set of isomorphism classes of irreducible graded
$K_\La$-modules.
For any graded $K_\La$-module $M$ we let
${M}^\circledast$ denote its graded dual, which means that
$({M}^\circledast)_{j} := \hom_{\C}(M_{-j},\C)$ and $x \in K_\La$ acts on $f
\in {M}^\circledast$ by $(xf)(m) := f(x^* m)$, where $x^*$ denotes the image of
$x$ under (\ref{star}).
Clearly we have that \begin{equation}
{L(\la)}^\circledast
 \cong L(\la)
\end{equation}
for each $\la \in \La$.

For $\la \in \La$, let $P(\la) := K_\La e_\la$. This is a graded $K_\La$-module
with basis $$\left\{(\underline{\nu} \mu \overline\la)\:\big|\:\text{for all
}\nu,\mu \in \La \text{ such that }\nu \subset \mu \supset \la\right\}.$$ By
Lemma~\ref{notmany}, $P(\la)$ is locally finite dimensional and bounded below,
and it is globally finite dimensional if and only if $\defect(\La) < \infty$.
So $P(\la)$ belongs to the category $\mod{K_\La}$. In fact, it is the
projective cover of $L(\la)$ in this category, and the modules
\begin{equation}\label{pims}
\{P(\la)\langle j \rangle\:|\:\la \in \La, j\in \Z\}
\end{equation}
give a full set of projective indecomposable modules in $\mod{K_\La}$.

\phantomsubsection{Grothendieck groups and Cartan matrix}
Let $\rep{K_\La}$ and $\proj{K_\La}$ be the full subcategories of
$\mod{K_\La}$ consisting of finite dimensional modules and finitely generated
projective modules, respectively.
The Grothendieck groups $[\rep{K_\La}]$ and
$[\proj{K_\la}]$ are the free $\Z$-modules with bases given by the isomorphism
classes of the irreducible modules from (\ref{simples}) and the projective
indecomposable modules from (\ref{pims}), respectively. If we view all our
Grothendieck groups instead as $\Z[q,q^{-1}]$-modules so that
$$
q^j [M] := [M\langle j
\rangle],
$$
then $[\rep{K_\La}]$ and $[\proj{K_\La}]$
become the free $\Z[q,q^{-1}]$-modules with bases
$\{[L(\la)] \:|\:\la \in \La\}$
and
$\{[P(\la)]\:|\:\la \in \La\}$, respectively. The Grothendieck
group $[\mod{K_\La}]$ of the category $\mod{K_\La}$ itself is the completion of
$[\rep{K_\La}]$ to a $\Z((q))$-module, where $\Z((q))$ denotes the ring of
formal Laurent series in $q$ (the localisation of $\Z[[q]]$ at $q$).

The inclusions of $\rep{K_\La}$ and $\proj{K_\La}$ as subcategories
of $\mod{K_\La}$ induce $\Z[q,q^{-1}]$-module homomorphisms
$$
[\rep{K_\La}] \;\xhookrightarrow{i}\;[\mod{K_\La}]
\;\xhookleftarrow{j}\;[\proj{K_\La}].
$$
The first map $i$ is obviously injective.
To see that the second map $j$ is injective too,
introduce the {\em $q$-Cartan matrix}
\begin{equation}\label{cartanmat}
C_\La(q) = (c_{\la,\mu}(q))_{\la,\mu \in\La}
\end{equation}
where
\begin{equation}\label{cmat}
c_{\la,\mu}(q) := \sum_{j \in \Z} q^j \dim \hom_{K_\La}(P(\la),P(\mu))_j \in \Z((q)).
\end{equation}
In other words, we have that
\begin{equation}\label{decm}
[P(\mu)] = \sum_{\la \in \La} c_{\la,\mu}(q) [L(\la)],
\end{equation}
equality written in $[\mod{K_\La}]$.
Since
$$
\hom_{K_\La}(P(\la), P(\mu)) = \hom_{K_\La} (K_\La e_\la, K_\La e_\mu)= e_\la
K_\La e_\mu
$$
and $e_\la K_\La e_\mu$ has basis
$\left\{(\underline{\la}\nu\overline{\mu})\:\big|\:
\nu \in \La \text{ such that }\la \subset\nu\supset\mu\right\}$,
we have that
\begin{equation}\label{exp}
c_{\la,\mu}(q) =
\sum_{\la \subset \nu \supset \mu}
q^{\deg(\underline{\la}\nu\overline{\mu})}.
\end{equation}
Recalling Lemma~\ref{notmany}, it follows that $c_{\la,\mu}(q)$ is in fact a
polynomial in $q$ and its constant coefficient is $1$ if $\la = \mu$, and $0$
otherwise. So the Cartan matrix $C_\La(q)$ is invertible in the ring of $\La
\times \La$ matrices with entries in $\Z((q))$. This implies that the map $j$
is injective as claimed, indeed, it shows that $[\mod{K_\La}]$ can also be
described as the completion of $[\proj{K_\La}]$. If $|\La| < \infty$ then we
even have that $[\rep{K_\La}] = [\proj{K_\La}]$ as submodules of
$[\mod{K_\La}]$.

\phantomsubsection{\boldmath Cell modules and $q$-decomposition numbers}
Now we introduce {\em cell modules} following \cite{GL}. For $\mu
\in \La$, define $V(\mu)$ to be the vector space on basis
\begin{equation}
\left\{
(\underline{\la} \mu| \:\big|\: \text{for all }
\la \in \La\text{ such that }\la \subset \mu\right\}.
\end{equation}
Equivalently, by Lemma~\ref{l2}(i), this basis is the set
$$
\left\{
(c \mu | \:\big|\:
\text{for all oriented cup diagrams $c \mu$}\right\}.
$$
We put a grading on $V(\mu)$ by defining
the degree of the vector $(c \mu |$ to be $\deg(c\mu)$ as in \eqref{dp}.
Note by Lemma~\ref{notmany} that $V(\mu)$ is finite dimensional
if and only if $\defect(\La) < \infty$,
though it is always locally finite dimensional and bounded below.
We make $V(\mu)$ into a graded $K_\La$-module by declaring
for any basis vector $(a \la b)$ of $K_\La$ that
\begin{equation}\label{Actby}
(a \la b) (c \mu| :=
\left\{
\begin{array}{ll}
s_{a \la b}(\mu)  (a \mu|
&\text{if $b^* = c$ and $a \mu$ is oriented,}\\
0&\text{otherwise.}
\end{array}\right.
\end{equation}
where $s_{a \la b}(\mu) \in \{0,1\}$ is the scalar
defined by Theorem~\ref{cell2}.
We remark that the proofs of Theorems~\ref{cell}
and \ref{cell2} give an explicit algorithm to compute these
scalars in practise.
We have not yet checked that (\ref{Actby}) makes
$V(\mu)$ into a well-defined module, but that
will be explained in the proof of the
next theorem.

\begin{Theorem}
\label{qh1}
For $\la \in \La$,
enumerate the $2^{\defect(\la)}$ distinct elements of the set
$\{\mu \in \La\:|\:\mu \supset \la\}$
as $\mu_1,\mu_2,\dots,\mu_n = \la$
so that $\mu_i > \mu_j$ implies $i < j$.
Let $M(0) := \{0\}$ and for $i=1,\dots,n$ define
$M(i)$ to be the subspace of $P(\la)$
generated by $M(i-1)$ and the vectors
$$
\left\{
(c \mu_i \overline{\la} ) \:\big|\:
\text{for all oriented cup diagrams $c \mu_i$}\right\}.
$$
Then
$$
\{0\} = M(0) \subset M(1) \subset\cdots\subset M(n) = P(\la)
$$
is a filtration of $P(\la)$ as a $K_\La$-module such that
$$
M(i) / M(i-1) \cong V(\mu_i) \langle \deg(\mu_i
\overline{\la})\rangle
$$
for each $i=1,\dots,n$.
\end{Theorem}

\begin{proof}
Corollary~\ref{ideal2} implies that each $M(i)$
is a $K_\La$-submodule of $P(\la)$.
Since $P(\la)$ has basis
$\left\{(c \mu \overline{\la})\:\big|\:\text{for all oriented circle diagrams
$c \mu \overline{\la}$}\right\}$,
the map
$$
f_i:V(\mu_i)\langle \deg(\mu_i \overline{\la}) \rangle
\rightarrow M(i) / M(i-1),\qquad
(c \mu_i| \mapsto (c \mu_i \overline{\la}) + M(i-1)
$$
is a vector space isomorphism, and it is homogeneous of degree $0$
thanks to (\ref{dp}).
Recalling Theorem~\ref{cell2}, we use this isomorphism
to transport the $K_\La$-module structure
on $M(i) / M(i-1)$ to the vector space $V(\mu_i)$, to make
$V(\mu_i)$ into a well-defined $K_\La$-module such that
$$
(a \la b) (c \mu_i| :=
\left\{
\begin{array}{ll}
s_{a \la b}(\mu_i) (a \mu_i|
&\text{if $b^* = c$ and $a \mu_i$ is oriented,}\\
0&\text{otherwise.}
\end{array}\right.
$$
This is the same as the formula (\ref{Actby}) with $\mu = \mu_i$,
so the map $f_i$ is a $K_\La$-module isomorphism.
Since any $\mu \in \La$ arises as $\mu= \mu_i$ for some $\la$ and
some $i$
(e.g. one can take $\la = \mu$ and $i=n$),
we have proved at the same time that $V(\mu)$ itself is a
well-defined $K_\La$-module.
\end{proof}

\begin{Theorem}\label{qh2}
For $\mu \in \La$, let $N(j)$ be the submodule of $V(\mu)$ spanned by all
graded pieces of degree $\geq j$. Then
$$
V(\mu) = N(0) \supseteq N(1) \supseteq N(2) \supseteq \cdots
$$
is a
filtration of $V(\mu)$ as a $K_\La$-module such that
$$
N(j) / N(j+1) \cong \bigoplus_{\substack{\la  \subset \mu\text{\,with}\\
\deg(\underline{\la} \mu) = j}}
  L(\la) \langle j \rangle
$$
for each $j \geq 0$. Moreover, we have that $N(j) = 0$ for $j \gg 0$, i.e.
$V(\mu)$ is finite dimensional,  if and only if $\defect(\La) < \infty$.
\end{Theorem}

\begin{proof}
Since $K_\La$ is positively graded, it is immediate that each $N(j)$ is a
submodule. The quotient $N(j) / N(j+1)$ has basis
$$
\left\{(\underline\la \mu | + N(j+1) \:\big|\: \la \in \La\text{ such that }
\la \subset \mu\text{ and }\deg(\underline{\la} \mu) = j\right\}.
$$
We claim that the subspace of $N(j) / N(j+1)$ spanned by $(\underline\la \mu |
+ N(j+1)$ is isomorphic to $L(\la) \langle j \rangle$. To see that, we just
need to observe for any basis vector $(a \nu b)$ of $K_\La$ and $\la \subset
\mu$ that
$$
(a \nu b)
(\underline\la \mu | \equiv \left\{
\begin{array}{ll}
(\underline\la \mu |&\text{if $a \nu b = \underline \la \la \overline\la$,}\\
0&\text{otherwise},
\end{array}
\right.
$$
working modulo $N(j+1)$. This is clear if $\deg(a \nu b) > 0$ or if $b \neq
\overline{\la}$ by (\ref{Actby}). So we may assume that $a \nu b =
\underline{\la} \la \overline{\la}$ and need to show that $s_{\underline{\la}
\la \overline{\la}}(\mu) = 1$. This follows from Theorem~\ref{cell2}(iii). The
final statement of the theorem is immediate from Lemma~\ref{notmany}(ii).
\end{proof}

For $\la,\mu \in \La$, we define
\begin{equation}\label{dm}
d_{\la,\mu}(q)
:=
\left\{
\begin{array}{ll}
q^{\deg(\underline{\la} \mu)} &\text{if $\la \subset \mu$,}\\
0&\text{otherwise.}
\end{array}\right.
\end{equation}
We call the resulting matrix
\begin{equation}\label{decmat}
D_\La(q) = (d_{\la,\mu}(q))_{\la,\mu \in \La}
\end{equation}
the {\em $q$-decomposition matrix} because by
Theorems~\ref{qh1} and \ref{qh2} we have that
\begin{align}
[V(\mu)] &= \sum_{\la \in \La} d_{\la,\mu}(q) [L(\la)],\label{dmat}\\
[P(\la)] &= \sum_{\mu \in \La} d_{\la,\mu}(q) [V(\mu)],
\end{align}
both equalities written in the Grothendieck group $[\mod{K_\La}]$.
Recall from Lemma~\ref{l2} that $\la \subset \mu$ implies $\la \leq \mu$.
So the $q$-decomposition matrix $D_\La(q)$ is upper
unitriangular when rows and columns are ordered in some way refining the Bruhat order.
Here is an example:
\begin{equation}\label{qdec}
\begin{array}{c|cccccc}
&
\scriptstyle{\down\down\up\up}&
\scriptstyle{\down\up\down\up}&
\scriptstyle{\up\down\down\up}&
\scriptstyle{\down\up\up\down}&
\scriptstyle{\up\down\up\down}&
\scriptstyle{\up\up\down\down}\\
\hline
\scriptstyle{\down\down\up\up}&1&q&0&0&q&q^2\\
\scriptstyle{\down\up\down\up}&0&1&q&q&q^2&0\\
\scriptstyle{\up\down\down\up}&0&0&1&0&q&0\\
\scriptstyle{\down\up\up\down}&0&0&0&1&q&0\\
\scriptstyle{\up\down\up\down}&0&0&0&0&1&q\\
\scriptstyle{\up\up\down\down}&0&0&0&0&0&1
\end{array}
\end{equation}
On comparing (\ref{dm}) with (\ref{exp}) and using (\ref{dp}),
we get the following familiar factorisation of the Cartan matrix
$C_\La(q)$:
\begin{equation}\label{brauer}
C_\La(q) = D_\La(q) D_\La(q)^T.
\end{equation}

\phantomsubsection{Highest weight categories}
We are ready to prove
the following important result.
Note here it is essential to exclude
$\defect(\La) = \infty$, since in that case
the cell modules are not Artinian.

\begin{Theorem}\label{ghw}
Assuming $\defect(\La) < \infty$, the category $\rep{K_\La}$ is a positively
graded highest weight category with duality in the sense of Cline, Parshall and
Scott. For $\la \in \La$ (the weight poset), the irreducible, standard, costandard,
projective indecomposable and injective indecomposable objects are the modules
$L(\la), V(\la), {V(\la)}^\circledast$, $P(\la)$ and ${P(\la)}^\circledast$,
respectively.
\end{Theorem}

Before proving the theorem, we give the definition of a positively graded highest weight category
with duality following \cite[Definition 3.1]{CPS1}, \cite[(1.1)]{CPS2} and
\cite[(1.2)]{CPS3} as follows. This means a $\Z$-category $\mathcal C$ in the
sense of \cite[Appendix E]{AJS} with shift operators denoted $M \mapsto
M\langle j \rangle$ for each $j \in \Z$ (like the category of finite
dimensional graded modules over a graded algebra) together with an
interval-finite poset $\La$ such that
\begin{itemize}
\item[(1)]
$\mathcal C$ is an $\C$-linear Artinian category equipped
with a duality, that is, an equivalence of categories $\circledast:\mathcal C \rightarrow
\mathcal C^{\op}$, such that $({M \langle j \rangle})^\circledast \cong
{M}^\circledast \langle -j \rangle$ for $M \in \mathcal C$ and $j \in \Z$;
\item[(2)]
for each $\la \in \La$ there is a given object $L(\la)
\cong{L(\la)}^\circledast \in\mathcal C$ such that $\{L(\la) \langle j
\rangle\:|\:\la \in \La, j \in \Z\}$ is a complete set of representatives for
the isomorphism classes of irreducible objects in $\mathcal C$;
\item[(3)]
each $L(\la)$ has a projective cover $P(\la)\in\mathcal C$ such that
all composition factors of $P(\la)$ are isomorphic to $L(\mu)\langle j \rangle$'s
for $\mu \in \La$ and $j \geq 0$;
\item[(4)] writing $V(\la)$ for the largest quotient of
$P(\la)$ with the property that all composition factors of
its radical are isomorphic to
$L(\mu)\langle j \rangle$'s for $\mu < \la$,
the object $P(\la)$ has a
filtration with $V(\la)$ at the top and all
other factors isomorphic to $V(\nu)\langle j \rangle$'s for
$\nu > \la$.
\end{itemize}

\begin{proof}[Proof of Theorem~\ref{ghw}]
We have to verify the above properties for $\rep{K_\La}$.
Property (1) is clear since
all the objects are finite dimensional modules. We have already noted that the
modules $L(\la)\langle j \rangle$ for $\la \in \La$ and $j \in \Z$ give a
complete set of irreducible graded $K_\La$-modules giving (2). For (3), we have
already noted that $P(\la)$ is the projective cover of $L(\la)$, so we are done
by (\ref{decm}) and (\ref{exp}). Finally, for (4), Theorem~\ref{qh1} shows that
the projective cover $P(\la)$ of $L(\la)$ has a filtration with top quotient
isomorphic to the cell module $V(\la)$ and all other quotients isomorphic to
$V(\mu)\langle j \rangle$'s for $\mu > \la$ and $j \in \Z$. Note in particular
this means that $V(\la)$ has a unique irreducible quotient, isomorphic to
$L(\la)$. It just remains to check that the cell module $V(\la)$ defined
earlier is isomorphic to the module $V(\lambda)$ defined in (4). In the notation of
Theorem~\ref{qh1}, $V(\la)$ is isomorphic to the quotient $Q := P(\la) / M(n-1)$,
and Theorem~\ref{qh2} shows that all composition factors of the radical of $Q$
are isomorphic to $L(\mu)\langle j \rangle$'s for $\mu < \la$. Moreover, $M(n-1)$
has a filtration with sections isomorphic to $V(\mu) \langle j \rangle$'s for
$\mu > \la$, so all irreducible quotients of $M(n-1)$ are isomorphic to
$L(\mu)\langle j \rangle$'s for $\mu > \la$. This completes the proof.
\end{proof}

If $|\La| < \infty$, so that $K_\La$ is a finite dimensional
algebra,
Theorem~\ref{ghw} is well known to be equivalent to the following statement;
see \cite[Theorem 3.6]{CPS1} and \cite[Proposition E.4]{AJS}.

\begin{Corollary}\label{gqh}
Assuming $|\La| < \infty$, the algebra $K_\La$ is a
positively graded quasi-hereditary algebra in the sense of Cline, Parshall and Scott.
\end{Corollary}

\section{Further results}
\label{section6}
In this section, we begin by explaining an effective direct way
to compute the multiplication in $K_\La$, bypassing the
closure and extension constructions from $\S$\ref{sstroppel}.
Then we generalise the definition of Khovanov's
algebra slightly, introducing an analogous cellular algebra $H_\La$
for any block $\La$ (not just Khovanov blocks).
Finally we show that $H_\La$ is a symmetric algebra,
extending a known result for Khovanov blocks,
and discuss briefly some other aspects of the representation theory
of generalised Khovanov algebras.

\phantomsubsection{The generalised surgery procedure}
We now introduce a generalisation of the surgery
procedure from $\S$\ref{skhovanov}
which makes sense for arbitrary (not necessarily closed) circle diagrams.
Suppose we are given a consistently oriented
diagram consisting of an oriented
cup diagram at the bottom,
a symmetric middle section, and then an oriented cap diagram
at the top.
Pick a symmetric pair of a cup and a cap in the middle section
that can be connected without crossings.
Cut open the cup and cap and stitch the loose ends
together exactly as in the original surgery procedure.
If the cup and cap both belonged to circles,
we re-orient the resulting
diagram(s)
following exactly the same rules
(\ref{mult})--(\ref{comult})
as before.
However it is now also possible that one or both of the cup and
cap belonged to a line, not a circle, before the cut was made.
In that case, we supplement the rules
(\ref{mult})--(\ref{comult})
with the following extra rules, in which
$y$ denotes a line and $1$ and $x$
denote anti-clockwise and clockwise circles as before:
\begin{align}
1 \otimes y &\mapsto y,\quad
x \otimes y \mapsto 0,\quad
y \mapsto x \otimes y,\\
y \otimes y &\mapsto
\left\{
\begin{array}{ll}
y\otimes y&\text{if both rays from one of the lines are oriented $\up$}\\
&\text{and both rays from the other line are oriented $\down$;}\\
0&\text{otherwise.}\label{yuk}
\end{array}
\right.
\end{align}
So, the rule $1 \otimes y\mapsto y$ indicates that an anti-clockwise circle
and a line convert to a single line,
the rule $y \mapsto x \otimes y$ indicates that a single line converts
to a clockwise circle and a line,
the rule $y \otimes y \mapsto y \otimes y$ indicates that
two lines convert to two (different) lines,
and so on.
It is important to note that
there is only ever one admissible
choice for the orientation
of vertices lying on lines; in particular,
the orientations of vertices at the ends of rays
are never changed by the generalised surgery procedure.

\begin{Theorem}\label{gs}
For any block $\La$, the product
$(a\la b) (b^* \mu d)$
of basis
vectors from $K_\La$ can be computed as follows.
\begin{itemize}
\item[(i)]
First draw $a \la b$ underneath $b^* \mu d$
and stitch corresponding rays together to create a new diagram
with a symmetric middle section;
corresponding rays in $\la b$ and $b^* \mu$
are necessarily oriented in the same direction so this makes sense.
\item[(ii)]
Then iterate the generalised surgery procedure just described
in any order that makes sense
to create a disjoint union of diagrams, none of which
has any cup-cap pairs left in its symmetric middle section.
\item[(iii)]
Finally
identify the top and bottom number lines of each
of the resulting diagrams to obtain a disjoint union of some new
oriented circle diagrams.
The desired product $(a \la b) (b^* \mu d)$ is the sum of the
corresponding basis vectors of $K_\La$.
\end{itemize}
\end{Theorem}

Before proving the theorem, we illustrate the statement by using it
to compute the following product:
\vspace{4.5mm}
\begin{equation}
\phantom{x}\label{step1}
\end{equation}
$$
\begin{picture}(265,4)
\put(0,20){\line(1,0){115}}
\put(128,17){$\times$}
\put(150,20){\line(1,0){115}}

\put(20.1,20.1){$\scriptstyle\down$}
\put(66.1,20.1){$\scriptstyle\down$}
\put(112.1,20.1){$\scriptstyle\down$}
\put(43.1,15.5){$\scriptstyle\up$}
\put(89.1,15.5){$\scriptstyle\up$}
\put(-2.9,15.5){$\scriptstyle\up$}

\put(80.5,20){\oval(23,23)[b]}
\put(80.5,20){\oval(23,23)[t]}
\put(34.5,20){\oval(23,23)[t]}
\put(80.5,20){\oval(69,40)[b]}
\put(0,0){\line(0,1){30}}
\put(23,0){\line(0,1){20}}
\put(115,20){\line(0,1){10}}

\put(161.5,20){\oval(23,23)[t]}
\put(184.5,20){\oval(23,23)[b]}
\put(230.5,20){\oval(23,23)[b]}
\put(150,0){\line(0,1){20}}
\put(265,0){\line(0,1){30}}
\put(242,20){\line(0,1){10}}
\put(219,20){\line(0,1){10}}
\put(196,20){\line(0,1){10}}

\put(170.1,20.1){$\scriptstyle\down$}
\put(239.1,20.1){$\scriptstyle\down$}
\put(262.1,20.1){$\scriptstyle\down$}
\put(193.1,15.5){$\scriptstyle\up$}
\put(216.1,15.5){$\scriptstyle\up$}
\put(147.1,15.5){$\scriptstyle\up$}
\end{picture}
$$
Putting the first diagram under the second
and stitching corresponding rays together
according to (i) gives:

\vspace{6mm}
\begin{equation}
\phantom{x}\label{step2}
\end{equation}
$$
\begin{picture}(120,28)
\put(0,20){\line(1,0){115}}
\put(0,55){\line(1,0){115}}

\put(20.1,20.1){$\scriptstyle\down$}
\put(66.1,20.1){$\scriptstyle\down$}
\put(112.1,20.1){$\scriptstyle\down$}
\put(43.1,15.5){$\scriptstyle\up$}
\put(89.1,15.5){$\scriptstyle\up$}
\put(-2.9,15.5){$\scriptstyle\up$}

\put(80.5,20){\oval(23,23)[b]}
\put(80.5,20){\oval(23,23)[t]}
\put(34.5,20){\oval(23,23)[t]}
\put(80.5,20){\oval(69,40)[b]}
\put(0,0){\line(0,1){55}}
\put(23,0){\line(0,1){20}}
\put(115,20){\line(0,1){45}}

\put(11.5,55){\oval(23,23)[t]}
\put(34.5,55){\oval(23,23)[b]}
\put(80.5,55){\oval(23,23)[b]}
\put(92,55){\line(0,1){10}}
\put(69,55){\line(0,1){10}}
\put(46,55){\line(0,1){10}}

\put(20.1,55.1){$\scriptstyle\down$}
\put(89.1,55.1){$\scriptstyle\down$}
\put(112.1,55.1){$\scriptstyle\down$}
\put(43.1,50.5){$\scriptstyle\up$}
\put(-2.9,50.5){$\scriptstyle\up$}
\put(66.1,50.5){$\scriptstyle\up$}
\end{picture}
$$
Then we make two iterations of the generalised surgery procedure
according to (ii):
$$
\begin{picture}(260,66)
\put(126,36){$\rightsquigarrow$}

\put(0,20){\line(1,0){115}}
\put(0,55){\line(1,0){115}}
\put(20.1,20.1){$\scriptstyle\down$}
\put(66.1,20.1){$\scriptstyle\down$}
\put(112.1,20.1){$\scriptstyle\down$}
\put(43.1,15.5){$\scriptstyle\up$}
\put(89.1,15.5){$\scriptstyle\up$}
\put(-2.9,15.5){$\scriptstyle\up$}
\put(80.5,20){\oval(23,23)[b]}
\put(80.5,20){\oval(23,23)[t]}
\put(34.5,20){\oval(23,23)[t]}
\put(80.5,20){\oval(69,40)[b]}
\put(0,0){\line(0,1){55}}
\put(23,0){\line(0,1){20}}
\put(115,20){\line(0,1){45}}
\put(11.5,55){\oval(23,23)[t]}
\put(34.5,55){\oval(23,23)[b]}
\put(80.5,55){\oval(23,23)[b]}
\put(92,55){\line(0,1){10}}
\put(69,55){\line(0,1){10}}
\put(46,55){\line(0,1){10}}
\put(20.1,55.1){$\scriptstyle\down$}
\put(89.1,55.1){$\scriptstyle\down$}
\put(112.1,55.1){$\scriptstyle\down$}
\put(43.1,50.5){$\scriptstyle\up$}
\put(-2.9,50.5){$\scriptstyle\up$}
\put(66.1,50.5){$\scriptstyle\up$}

\put(150,20){\line(1,0){115}}
\put(150,55){\line(1,0){115}}
\put(170.1,20.1){$\scriptstyle\down$}
\put(239.1,20.1){$\scriptstyle\down$}
\put(262.1,20.1){$\scriptstyle\down$}
\put(193.1,15.5){$\scriptstyle\up$}
\put(216.1,15.5){$\scriptstyle\up$}
\put(147.1,15.5){$\scriptstyle\up$}
\put(230.5,20){\oval(23,23)[b]}
\put(230.5,20){\oval(69,40)[b]}
\put(150,0){\line(0,1){55}}
\put(173,0){\line(0,1){55}}
\put(265,20){\line(0,1){45}}
\put(161.5,55){\oval(23,23)[t]}
\put(242,20){\line(0,1){45}}
\put(219,20){\line(0,1){45}}
\put(196,20){\line(0,1){45}}
\put(170.1,55.1){$\scriptstyle\down$}
\put(239.1,55.1){$\scriptstyle\down$}
\put(262.1,55.1){$\scriptstyle\down$}
\put(193.1,50.5){$\scriptstyle\up$}
\put(147.1,50.5){$\scriptstyle\up$}
\put(216.1,50.5){$\scriptstyle\up$}

\dashline{2}(80.5,32)(80.5,43)
\dashline{2}(34.5,32)(34.5,43)
\end{picture}
$$
Finally we identify top and bottom number lines as in
(iii) to get the answer:
$$
\begin{picture}(125,34)
\put(0,20){\line(1,0){115}}
\put(11.5,20){\oval(23,23)[t]}
\put(80.5,20){\oval(23,23)[b]}
\put(80.5,20){\oval(69,40)[b]}
\put(23,0){\line(0,1){20}}
\put(0,0){\line(0,1){20}}
\put(115,20){\line(0,1){10}}
\put(92,20){\line(0,1){10}}

\put(69,20){\line(0,1){10}}
\put(46,20){\line(0,1){10}}
\put(20.1,20.1){$\scriptstyle\down$}
\put(89.1,20.1){$\scriptstyle\down$}
\put(112.1,20.1){$\scriptstyle\down$}
\put(43.1,15.5){$\scriptstyle\up$}
\put(66.1,15.5){$\scriptstyle\up$}
\put(-2.9,15.5){$\scriptstyle\up$}
\end{picture}
$$
We leave it to the reader to check that this is
indeed the desired product according to the
original definition of multiplication.

\begin{proof}[Proof of Theorem~\ref{gs}]
It is enough to consider the case that
$\La$ consists of bounded weights.
We require an obvious extension of the notion of the closure
 of a diagram from (\ref{notion})
to diagrams
with a symmetric middle section. For example, the closure
of the diagram (\ref{step2}) is:
\vspace{16mm}
\begin{equation}
\phantom{x}\label{step3}
\end{equation}
$$
\begin{picture}(120,44)
\put(-46,40){\line(1,0){208}}
\put(-46,75){\line(1,0){208}}

\put(20.1,40.1){$\scriptstyle\down$}
\put(66.1,40.1){$\scriptstyle\down$}
\put(112.1,40.1){$\scriptstyle\down$}
\put(43.1,35.5){$\scriptstyle\up$}
\put(89.1,35.5){$\scriptstyle\up$}
\put(-2.9,35.5){$\scriptstyle\up$}

%\put(-25.9,40.1){$\scriptstyle\down$}
%\put(-48.9,40.1){$\scriptstyle\down$}
%\put(-25.9,75.1){$\scriptstyle\down$}
%\put(-48.9,75.1){$\scriptstyle\down$}

%\put(135.1,35.5){$\scriptstyle\up$}
%\put(158.1,35.5){$\scriptstyle\up$}
%\put(135.1,70.5){$\scriptstyle\up$}
%\put(158.1,70.5){$\scriptstyle\up$}

\put(11.5,75){\oval(69,40)[t]}
\put(11.5,75){\oval(115,57)[t]}
\put(126.5,75){\oval(69,40)[t]}
\put(126.5,75){\oval(23,23)[t]}

\put(-11.5,40){\oval(23,23)[b]}
\put(80.5,40){\oval(115,57)[b]}
\put(57.5,40){\oval(207,74)[b]}

\put(80.5,40){\oval(23,23)[b]}
\put(80.5,40){\oval(23,23)[t]}
\put(34.5,40){\oval(23,23)[t]}
\put(80.5,40){\oval(69,40)[b]}
\put(0,40){\line(0,1){35}}
\put(115,40){\line(0,1){35}}

\put(138,40){\line(0,1){35}}
\put(161,40){\line(0,1){35}}
\put(-23,40){\line(0,1){35}}
\put(-46,40){\line(0,1){35}}

\put(11.5,75){\oval(23,23)[t]}
\put(34.5,75){\oval(23,23)[b]}
\put(80.5,75){\oval(23,23)[b]}

\put(20.1,75.1){$\scriptstyle\down$}
\put(89.1,75.1){$\scriptstyle\down$}
\put(112.1,75.1){$\scriptstyle\down$}
\put(43.1,70.5){$\scriptstyle\up$}
\put(-2.9,70.5){$\scriptstyle\up$}
\put(66.1,70.5){$\scriptstyle\up$}
\end{picture}
$$
We will refer to the vertices in the diagram before
closing as {\em inner vertices} and the additional vertices
that get added as {\em outer vertices}.
In (\ref{step3}), we have omitted the $\down$'s labelling
outer vertices on the left hand side and the $\up$'s labelling
outer vertices on the right hand side; this
convention helps to
distinguish inner and outer vertices in practise.

According to the definition of multiplication in $K_\La$,
to compute
the product $(a \la b) (b^* \mu d)$ we should
first close both diagrams, then
draw one under the other and iterate the surgery procedure.
Finally at the end we identify top and bottom number lines
as usual, discarding
any diagrams corresponding to basis vectors in the
ideal $I_\La$, then remove the outer vertices in the remaining diagrams
to get the final answer.
If at any point during the sequence of surgery procedures
some diagram is produced with an $\up$ labeling
an outer vertex on the left or a $\down$ labeling an outer vertex
on the right, we can simply omit
that diagram from the outset, since it can only ever create
diagrams belonging to the ideal $I_\La$ at the end.
There are three
types of surgery procedure needed in this process:
\begin{itemize}
\item[(a)] surgeries involving
caps which pass though two outer vertices;
\item[(b)] surgeries involving
caps which pass though one
outer vertex;
\item[(c)] surgeries involving
caps which pass though no outer vertices.
\end{itemize}
We order the sequence of surgery
procedures so that
all the ones of type (a) are performed first, followed by
type (b), followed by type (c).

The surgeries of type (a) are rather trivial: they only ever involve
the rule $1 \otimes 1 \mapsto 1$ and no vertices change orientation
during these operations. Now assume all type (a) surgeries have
been performed and consider the surgeries of type (b).
For example if we are in the situation of (\ref{step1})
and have just performed the one required type (a) surgery, we obtain
the following diagram in which two type (b) surgeries are needed:
\vspace{16mm}
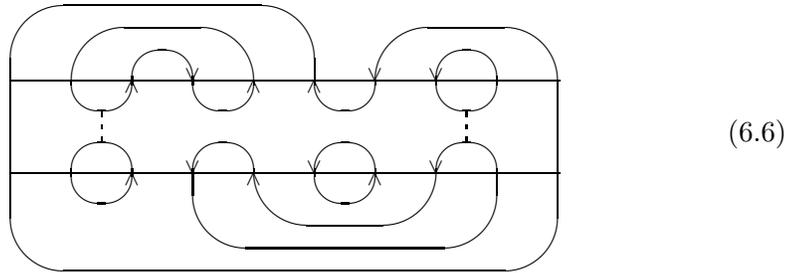
\begin{equation}
\phantom{x}\label{step4}
\end{equation}
$$
\begin{picture}(120,46)
\put(-46,40){\line(1,0){208}}
\put(-46,75){\line(1,0){208}}

\put(20.1,40.1){$\scriptstyle\down$}
\put(66.1,40.1){$\scriptstyle\down$}
\put(112.1,40.1){$\scriptstyle\down$}
\put(43.1,35.5){$\scriptstyle\up$}
\put(89.1,35.5){$\scriptstyle\up$}
\put(-2.9,35.5){$\scriptstyle\up$}

%\put(-25.9,40.1){$\scriptstyle\down$}
%\put(-48.9,40.1){$\scriptstyle\down$}
%\put(-25.9,75.1){$\scriptstyle\down$}
%\put(-48.9,75.1){$\scriptstyle\down$}

%\put(135.1,35.5){$\scriptstyle\up$}
%\put(158.1,35.5){$\scriptstyle\up$}
%\put(135.1,70.5){$\scriptstyle\up$}
%\put(158.1,70.5){$\scriptstyle\up$}

\put(11.5,75){\oval(69,40)[t]}
\put(11.5,75){\oval(115,57)[t]}
\put(126.5,75){\oval(69,40)[t]}
\put(126.5,75){\oval(23,23)[t]}

\put(-11.5,40){\oval(23,23)[b]}
\put(80.5,40){\oval(115,57)[b]}
\put(57.5,40){\oval(207,74)[b]}

\put(-11.5,40){\oval(23,23)[t]}
\put(-11.5,75){\oval(23,23)[b]}

\put(80.5,40){\oval(23,23)[b]}
\put(126.5,40){\oval(23,23)[t]}
\put(126.5,75){\oval(23,23)[b]}
\put(80.5,40){\oval(23,23)[t]}
\put(34.5,40){\oval(23,23)[t]}
\put(80.5,40){\oval(69,40)[b]}

\put(161,40){\line(0,1){35}}
\put(-46,40){\line(0,1){35}}

\put(11.5,75){\oval(23,23)[t]}
\put(34.5,75){\oval(23,23)[b]}
\put(80.5,75){\oval(23,23)[b]}

\put(20.1,75.1){$\scriptstyle\down$}
\put(89.1,75.1){$\scriptstyle\down$}
\put(112.1,75.1){$\scriptstyle\down$}
\put(43.1,70.5){$\scriptstyle\up$}
\put(-2.9,70.5){$\scriptstyle\up$}
\put(66.1,70.5){$\scriptstyle\up$}

\dashline{2}(126.5,52)(126.5,63)
\dashline{2}(-11.5,52)(-11.5,63)

\end{picture}
$$
The surgeries of type (b) are obviously
in bijective correspondence
with the pairs of rays that need to be stitched together
according to (i) from the statement of the theorem.
In fact,
after ignoring diagrams with an $\up$ labelling
an outer vertex on the left or a $\down$ labelling an outer vertex
on the right,
there is a unique diagram
obtained after performing all type
(b) surgeries and this diagram coincides with the closure
of the diagram obtained by stitching rays together
according to (i);
e.g. in the running example we get exactly the diagram
(\ref{step3}) on making the two type (b) surgeries
to (\ref{step4}).
This follows in general
by considering the types of lines
that the rays in $b$ and $b^*$ can possibly belong to, and considering
their closures in each case.

So now we have performed all surgeries of types (a) and (b),
and have obtained a diagram which is simply the closure
of the diagram obtained after performing (i).
Finally, we need to perform the surgeries of type (c).
These are in bijective correspondence
with the generalised surgeries required in (ii).
To complete the proof of the theorem,
we run the two sequences of surgeries (from (ii) and from (c)) in parallel.
The point is that the closures of the diagrams obtained
after each generalised surgery from (ii) are exactly the same
as the diagrams obtained after performing the corresponding
surgery of type (c)
 (ignoring diagrams with an $\up$ on the left or a $\down$
on the right).
This follows by considering the various
cases that occur in a single generalised surgery procedure
and comparing each with its closure. We omit the details.
\end{proof}

\phantomsubsection{\boldmath The generalised Khovanov algebra}
Let $\La$ be any block
and recall the definition (\ref{blockdefect}) of $\defect(\La)$.
Let
\begin{equation}
\La^\circ := \{\la \in \La\:|\:\defect(\la) = \defect(\La)\}
\end{equation}
denote the subset of $\La$ consisting of all weights of
{\em maximal defect}.
If $\La$ is a Khovanov block this is the same as the set $\La^\circ$ defined
at the beginning of $\S$\ref{skhovanov}.
Define the {\em generalised Khovanov algebra} $H_\La$ by setting
\begin{equation}
H_\La := \bigoplus_{\al, \be \in \La^\circ} e_\al \, K_\La\, e_\be.
\end{equation}
If $|\La| < \infty$, we have simply that
$H_\La = e K_\La e$ where
$e := \sum_{\la \in \La^\circ} e_\la \in K_\La$.
In general, $H_\La$ has basis
\begin{equation}
\{(\underline{\al} \la \overline{\be})\:|\:
\text{for all $\al \in \La^\circ, \la \in \La,
\be \in \La^\circ$  such that }\al \subset \la \supset \be\}.
\end{equation}
Moreover, the degree $0$ component of $H_\La$ has basis
\begin{equation}\label{ela3}
\{e_\la\:|\:\la \in \La^\circ\},
\end{equation}
recalling (\ref{ela2}).
The multiplication in $H_\La$
can be computed
easily using Theorem~\ref{gs} and the generalised surgery procedure;
moreover the
rule $y \otimes y \mapsto y \otimes y$
from (\ref{yuk}) is never needed in this maximal defect
situation.

If $\La$ is a Khovanov block it is straightforward to see that
$H_\La$ as defined here is the same as Khovanov's algebra
$H_\La$ from $\S$\ref{skhovanov}.
If $\defect(\La) = \infty$ then $H_\La = \{0\}$;
{\em we will exclude this boring
case for the remainder of the section.}
Finally if $|\La| = \infty$
and $\defect(\La) < \infty$
it often happens
that $H_\La = K_\La$. More precisely, using notation from the introduction,
it is the case that
\begin{equation}
H_\La = K_\La
\end{equation}
if $K_\La \cong
K^\infty_m \cong K^m_\infty$ or
$K_\La \cong K^{+\infty}_m \cong K^m_{-\infty}$,
but not if $K_\La \cong
K^{-\infty}_m \cong K^m_{+\infty}$.
(The isomorphisms
$K^\infty_m \cong K^m_\infty$,
$K^{+\infty}_m \cong K^m_{-\infty}$
and
$K^{-\infty}_m \cong K^m_{+\infty}$
noted here come from
\eqref{star} and \eqref{curve}.)

The results from $\S$\ref{skhovanov} extend
almost at once to the generalised Khovanov algebra $H_\La$
for any $\La$.
For example the anti-automorphism from
(\ref{star}) restricts to an anti-automorphism
\begin{align}\label{starnew}
*:H_\La &\rightarrow H_\La,
\qquad\quad
(a \la b) \mapsto (b^* \la a^*).
\end{align}
The following theorem
generalises Corollary~\ref{iscell}.

\begin{Theorem}\label{iscell3}
Assuming $|\La| < \infty$,
the algebra $H_\La$
is a cellular algebra in the sense of Graham and Lehrer with cell datum
$(\La, M, C, *)$ where
\begin{itemize}
\item[(i)]
$M(\la)$ denotes $\left\{\alpha \in \La^\circ\:|\:\alpha\subset \la\right\}$
for each $\la \in \La$;
\item[(ii)]
$C$ is defined by setting
$C^\la_{\alpha,\beta}:=( \underline{\alpha} \la \overline{\beta} )$
for $\la \in \La$ and $\alpha,\beta \in M(\la)$;
\item[(iii)]
$*$ is the anti-automorphism from (\ref{starnew}).
\end{itemize}
\end{Theorem}

\begin{proof}
Argue as in the proof of Corollary~\ref{iscell}
using Theorem~\ref{cell2}.
\end{proof}

\phantomsubsection{Symmetric algebra structure}
Now we prove that the generalised Khovanov algebra $H_\La$
is a symmetric algebra, In fact, we prove that it possesses a
{\em symmetrising form}
$\tau:H_\La \rightarrow \C$ such that the
induced bilinear form
$H_\La \otimes H_\La \rightarrow \C, a \otimes b \mapsto \tau(ab)$
is non-degenerate and symmetric.
For Khovanov blocks, this was proved
already by Khovanov in
\cite[Proposition 32]{K2}.

For a basis vector $(a \la b) \in H_\La$, we let
$(a \la b)^\#$ denote the basis vector of $H_\La$
whose diagram is obtained
by reversing the orientation of every circle in the diagram
$b^* \la a^*$.
The map $(a \la b) \mapsto (a \la b)^\#$
is an involution on the basis of $H_\La$
such that
$\deg(a \la b)^\# = 2 \defect(\La) - \deg (a \la b)$.
In particular, this shows that the top degree component
of $H_\La$ is in degree
$2 \defect(\La)$ and by (\ref{ela3}) it
has basis $\{e_\la^\#\:|\:\la \in \La^\circ\}$.

\begin{Theorem}\label{salg}
Let $\tau:H_\La \rightarrow \C$ be the linear map
such that $(a \la b) \mapsto 1$ if $\deg(a \la b) = 2 \defect(\La)$
and $(a \la b) \mapsto 0$ otherwise.
Then
$$
\tau((a \la b) (c \mu d)) = \left\{
\begin{array}{ll}
1&\text{if $(c \mu d) = (a \la b)^\#$,}\\
0&\text{otherwise.}
\end{array}\right.
$$
Hence, $\tau$ is a symmetrising form
and $H_\La$ is a symmetric algebra.
\end{Theorem}

\begin{proof}
We have that $(c \mu d) = (a \la b)^\#$
if and only if
\begin{itemize}
\item $d = a^*$ and $c = b^*$, i.e. the circle
diagrams $ab$ and $cd$ are mirror images of each other;
\item the mirror image pairs of circles in $ab$ and $cd$
are oriented in opposite ways in $a \la b$ and in $c \mu d$.
\end{itemize}
It is clear that $\tau((a \la b) (c \mu d)) = 0$ unless $d = a^*$ and $c =
b^*$. Assume therefore that $d = a^*$ and $c = b^*$.
Let $m$ be the number of circles in $a b$ and $n := \defect(\La)$.

The first surgery procedure performed to a mirror image pair of circles
in $ab$ and $cd$ turns two circles into one circles.
So after making exactly $m$ surgery procedures, one to each
mirror image pair of circles in $ab$ and $cd$,
there will be a total of $m$ circles, and still $(n-m)$
more surgery procedures to do.
Now $\tau((a \la b) (c \mu d))\not=0$ only
if there are exactly $n$ clockwise circles at the
end. To get $n$
circles at all, each of the remaining $(n-m)$ surgeries
have to turn one circle into two.
The final circles then get oriented via the rules $1 \mapsto 1
\otimes x + x \otimes 1$ or $x \mapsto x \otimes x$. It follows that to
get $n$ clockwise
circles at the end, the $m$ circles obtained after the first $m$ surgeries
have to be all of type $x$, i.e. each mirror image pair of circles
must indeed have been oriented in opposite ways in $a \la b$ and $c \mu d$
at the start.
\end{proof}

\phantomsubsection{Representations of the generalised Khovanov algebra}
To conclude the article, we indicate briefly how to deduce
the basic results about the representation theory of $H_\La$
from the analogous results already established for $K_\La$.
To do this, we exploit the exact functor
\begin{equation}\label{efunc}
e:\mod{K_\La} \rightarrow \mod{H_\La}
\end{equation}
mapping a module $M$ to
$eM := \bigoplus_{\la \in \La^\circ} e_\la M$
and defined by restriction on morphisms.
When $|\La| < \infty$ this is just the obvious
truncation functor arising from the
idempotent $e = \sum_{\la \in \La^\circ} e_\la$.
For $\la \in \La$,
we have that
$e L(\la) \neq 0$ if and only if $\la \in \La^\circ$,
and
$\{e L(\la)\:|\:\la \in \La^\circ\}$
is a complete set of irreducible $H_\La$-modules (up to degree shifts).
The modules
$\{e V(\mu)\:|\:\mu \in \La\}$
are the {\em cell modules} for $H_\La$
in the sense of Graham and Lehrer
arising from the cellular structure from Theorem~\ref{iscell3}.
Finally we refer to the modules
$
\{e P(\la)\:|\:\la \in \La\}
$
as {\em Young modules}, by analogy with the modular representation theory of
symmetric groups. If $\la \in \La^\circ$ then $e P(\la)$
is the projective cover of $e L(\la)$.
Applying the exact functor $e$,
the results on transition matrices in the graded Grothendieck group
from $\S$\ref{sdecomposition} easily imply similar statements
for $H_\La$.

\end{document}